\documentclass[aos,seceqn,citesort,dvips]{arximspdf}
\usepackage{graphics}
\doi{10.1214/09-AOS790}
\volume{38}
\issue{4}
\pubyear{2010}
\firstpage{2525}
\lastpage{2558}

\makeatletter

\newtheorem{theorem}{Theorem}
\newtheorem{corollary}{Corollary}
\newtheorem{lemma}{Lemma}
\newproclaim{remark}{Remark}
\newproclaim{Condition}{Condition}

\makeatother

\begin{document}
\begin{frontmatter}

\title{SPADES and mixture models\protect\thanksref{T1}}
\runtitle{SPADES and mixture models}

\thankstext{T1}{Part of the research was done while the authors were
visiting the Isaac
Newton Institute for Mathematical Sciences (Statistical Theory and
Methods for Complex, High-Dimensional Data Programme) at Cambridge
University during Spring 2008.}

\begin{aug}
\author[A]{\fnms{Florentina} \snm{Bunea}\thanksref{t2}\ead[label=e1]{bunea@stat.fsu.edu}},
\author[B]{\fnms{Alexandre B.} \snm{Tsybakov}\thanksref{t3}\ead[label=e2]{alexandre.tsybakov@upmc.fr}},\\
\author[A]{\fnms{Marten H.} \snm{Wegkamp}\corref{}\thanksref{t2}\ead[label=e3]{wegkamp@stat.fsu.edu}} and
\author[A]{\fnms{Adrian} \snm{Barbu}\ead[label=e4]{abarbu@stat.fsu.edu}}
\runauthor{Bunea, Tsybakov, Wegkamp and Barbu}
\affiliation{Florida State University, Laboratoire de Statistique, CREST and LPMA, Universit\'{e} Paris 6,
Florida State University and Florida State University}
\address[A]{F. Bunea\\
M. Wegkamp\\
A. Barbu\\
Department of Statistics\\
Florida State University\\
Tallahassee, Florida 32306-4330\\
USA\\
\printead{e1}\\
\phantom{E-mail: }\printead*{e3}\\
\phantom{E-mail: }\printead*{e4}} 
\address[B]{A. Tsybakov\\
Laboratoire de Statistique, CREST\\
92240 Malakoff\\
and\\
LPMA (UMR CNRS 7599)\\
Universit\'{e} Paris 6\\
75252 Paris, Cedex 05\\
France\\
\printead{e2}}
\end{aug}

\thankstext{t2}{Supported in part by the NSF Grant DMS-07-06829.}

\thankstext{t3}{Supported in part by the Grant ANR-06-BLAN-0194 and by
the PASCAL Network of
Excellence.}

\received{\smonth{1} \syear{2009}}
\revised{\smonth{12} \syear{2009}}

%
\begin{abstract}
This paper studies sparse density estimation via $\ell_1$
penalization (SPADES). We focus on estimation in high-dimensional
mixture models and nonparametric adaptive density estimation. We
show, respectively, that SPADES can recover, with high probability,
the unknown components of a mixture of probability densities and
that it yields minimax adaptive density estimates. These results are
based on a general sparsity oracle inequality
that the SPADES estimates satisfy. We offer a data driven method for
the choice of the tuning parameter used in the construction of SPADES.
The method uses the generalized bisection method first introduced in
\cite{bb09}. The suggested procedure bypasses the need for a grid
search and offers substantial computational savings. We complement our
theoretical results with a simulation study that employs this method
for approximations of one and two-dimensional densities with mixtures.
The numerical results strongly support our theoretical findings.
\end{abstract}

%
\begin{keyword}[class=AMS]
\kwd[Primary ]{62G08}
\kwd[; secondary ]{62C20}
\kwd{62G05}
\kwd{62G20}.
\end{keyword}
\begin{keyword}
\kwd{Adaptive estimation}
\kwd{aggregation}
\kwd{lasso}
\kwd{minimax risk}
\kwd{mixture models}
\kwd{consistent model selection}
\kwd{nonparametric
density estimation}
\kwd{oracle inequalities}
\kwd{penalized least squares}
\kwd{sparsity}
\kwd{statistical learning}.
\end{keyword}

\end{frontmatter}

\section{Introduction}

Let $X_1,\ldots,X_n$ be independent random variables with common
unknown density $f$ in $\mathbb R^d$. Let $\{f_1,\ldots,f_M\}$ be a
finite set of functions with $f_j\in L_2(\mathbb R^d), j=1,\ldots,M$,
called a dictionary. We consider estimators of $f$ that belong to
the linear span of $\{f_1,\ldots,f_M\}$. We will be particularly
interested in the case where $M \gg n$. Denote by $ \mathsf{f}_\lambda$
the linear combinations
\[
\mathsf{f}_\lambda(x) =\sum_{j=1}^M \lambda_j f_j(x),\qquad
\lambda=(\lambda_1,\ldots, \lambda_M) \in\mathbb R^M.
\]
Let us mention
some examples where such estimates are of importance:

\begin{itemize}
\item\textit{Estimation in sparse mixture models}.
Assume that the density $f$ can be represented as a finite mixture $
f=\mathsf{f}_{\lambda^*}$ where $f_j$ are known probability densities
and $\lambda^*$ is a vector of mixture probabilities. The number $M$
can be very large, much larger than the sample size $n$, but we
believe that the representation is sparse, that is, that very few
coordinates of $\lambda^*$ are nonzero, with indices corresponding to
a set $I^* \subseteq\{1, \ldots, M\}$. Our goal is to estimate the
weight vector
$\lambda^*$ by a vector ${\widehat\lambda}$ that adapts to this
unknown sparsity and
to identify $I^*$, with high probability.

\item\textit{Adaptive nonparametric density estimation.} Assume that the
density $f$ is a
smooth function, and $\{f_1,\ldots,f_M\}$ are the first $M$
functions from a basis in $L_2(\mathbb R^d)$. If the basis is orthonormal,
a natural idea is to estimate $f$ by an orthogonal series estimator
which has the form $\mathsf{f}_{\tilde\lambda}$ with $\tilde\lambda$
having the coordinates $\tilde\lambda_j=n^{-1}\sum_{i=1}^n
f_j(X_i)$. However, it is well known that such estimators are very
sensitive to the choice of $M$, and a data-driven selection of $M$
or thresholding is needed to achieve adaptivity (cf., e.g.,
\cite{rud82,kpt96,bm97}); moreover, these methods have been applied
with $M\le n$. We would like to cover more general problems where
the system $\{f_j\}$ is not necessarily orthonormal, even not
necessarily a basis, $M$ is not necessarily smaller than $n$, but an
estimate of the form $\mathsf{f}_{\widehat\lambda}$ still achieves,
adaptively, the optimal rates of convergence.

\item\textit{Aggregation of density estimators.} Assume now that
$f_1,\ldots,f_M$ are some
preliminary estimators of $f$ constructed from a training sample
independent of $(X_1,\ldots,X_n)$, and we would like to aggregate
$f_1,\ldots,f_M$. This means that we would like to construct a new
estimator, the aggregate, which is approximately as good as the best
among $f_1,\ldots,f_M$ or approximately as good as the best linear or
convex combination of $f_1,\ldots,f_M$. General notions of
aggregation and optimal rates are introduced in \cite{nem00,tsy03}.
Aggregation of density estimators is discussed in
\cite{st06,rt04,r06} and more recently in \cite{bir08} where one can
find further references. The aggregates that we have in mind here
are of the form $\mathsf{f}_{\widehat\lambda}$ with suitably chosen
weights ${\widehat\lambda}= {\widehat\lambda}(X_1,\ldots,X_n)\in
\mathbb R^M$.
\end{itemize}

In this paper we suggest a data-driven choice of
${\widehat\lambda}$ that can be used in all the examples mentioned
above and also more generally. We define ${\widehat\lambda}$ as a
minimizer of an $\ell_1$-penalized criterion, that we call SPADES
(SPArse density EStimation). This method was introduced in
\cite{us}. The idea of $\ell_1$-penalized estimation is widely used
in the statistical literature, mainly in linear regression where it
is usually referred to as the Lasso criterion
\cite{t96,cds01,det04,gr04,mb06}. For Gaussian sequence models or
for regression with an orthogonal design matrix the Lasso is equivalent
to soft thresholding \cite{d95,lg02}. Model selection consistency of
the Lasso type linear regression estimators is treated in many
papers including \cite{mb06,zy,zh,zou,l}. Recently,
$\ell_1$-penalized methods have been extended to nonparametric regression
with general fixed or random design \cite{btw05,btw06a,btw06b,brt},
as well as to some classification and other more general prediction
type models \cite{k05,k06,vdg06,bun08}.

In this paper we show that $\ell_1$-penalized techniques can also be
successfully used in density estimation. In Section \ref{sec2} we give the
construction of the SPADES estimates and we show that they satisfy
general oracle inequalities in Section \ref{sec3}. In the remainder of the
paper we discuss the implications of these results for two
particular problems, identification of mixture components and
adaptive nonparametric density estimation. For the application of
SPADES in aggregation problems we refer to \cite{us}.

Section \ref{sec4} is devoted to mixture models. A vast amount of literature
exists on estimation in mixture models, especially when the number
of components is known; see, for example, \cite{larry} for examples
involving the EM algorithm.
The literature on determining the number of mixture components
is still developing, and
we will focus on this aspect here. Recent works on the selection of
the number of components (mixture complexity) are
\cite{JPM01,BD05}. A consistent selection procedure specialized to Gaussian
mixtures is suggested in \cite{JPM01}. The method of \cite{JPM01}
relies on comparing a nonparametric kernel density estimator with
the best parametric fit of various given mixture complexities.
Nonparametric estimators based on the combinatorial density method
(see \cite{dl00}) are studied in \cite{BD05,BD08}. These can be
applied to estimating consistently the number of mixture components,
when the components have known functional form. Both
\cite{JPM01,BD05} can become computationally infeasible when $M$, the
number of
candidate components, is large. The method proposed here bridges
this gap and guarantees correct identification of the mixture
components with probability close to 1.

In Section \ref{sec4} we begin by giving conditions under which the mixture
weights can be estimated accurately, with probability close to 1.
This is an intermediate result that allows us to obtain the main
result of Section \ref{sec4}, correct identification of the mixture
components. We show that in identifiable mixture models, if the
mixture weights are above the noise level, then the components of
the mixture can be recovered with probability larger than $1 -
\varepsilon$, for any given small $\varepsilon$. Our results are
nonasymptotic, they hold for any $M$ and $n$. Since the emphasis
here is on correct component selection, rather than optimal density
estimation, the tuning sequence that accompanies the $\ell_1$
penalty needs to be slightly larger than the one used for good
prediction. The same phenomenon has been noted for
$\ell_1$-penalized estimation in linear and generalized regression models;
see, for example, \cite{bun08}.

Section \ref{sec5} uses the oracle inequalities of Section \ref{sec3} to show that
SPADES estimates adaptively achieve optimal rates of convergence (up
to a logarithmic factor) simultaneously on a large scale of
functional classes, such as H\"{o}lder, Sobolev or Besov classes, as
well as on the classes of sparse densities, that is, densities having
only a finite, but unknown, number of nonzero wavelet coefficients.

Section \ref{sec61} offers an algorithm for computing the SPADES. Our procedure
is based on
coordinate descent optimization, recently suggested by \cite{fhht07}.
In Section~\ref{sec:gbm} we use this algorithm together with a tuning parameter
chosen in a data adaptive manner.
This choice employs the generalized bisection method first introduced
in \cite{bb09}, a computationally efficient method for constructing
candidate tuning parameters without performing a grid search. The final
tuning parameter is chosen from the list of computed candidates by
using a 10-fold cross-validated dimension-regularized criterion.
The combined procedure works very well in practice, and we present a
simulation study in Section \ref{sec63}.

\section{Definition of SPADES}\label{sec2}
Consider the $L_2(\mathbb R^d)$ norm
\[
\Vert g \Vert= \biggl( \int_{\mathbb R^d} g^2(x) \,dx \biggr)^{1/2}
\]
associated with the inner product
\[
\langle g, h\rangle =\int_{\mathbb R^d} g(x) h(x) \,d x
\]
for $g,h\in L_2(\mathbb R^d)$. Note that if the density $f$ belongs to
$L_2(\mathbb R^d)$ and $X$ has the same distribution as $X_i$, we have,
for any $g\in L_2$,
\[
\langle g, f\rangle = \mathbb{E}g(X) ,
\]
where the expectation is taken under $f$. Moreover,
%
%
\begin{equation}\label{id1}
\| f- g \|^2 = \| f \|^2 + \|g\|^2 -2\langle g, f\rangle  = \|f\|^2 + \|g\|^2
-2\mathbb{E}g(X).
\end{equation}
In view of identity (\ref{id1}), minimizing $\| \mathsf{f}_\lambda- f
\|^2$ in $\lambda$ is the same as minimizing
\[
\gamma(\lambda) = -2 \mathbb{E}\mathsf{f}_\lambda(X)
+\|\mathsf{f}_\lambda\|^2.
\]
The function $\gamma(\lambda)$ depends on
$f$ but can be approximated by its empirical counterpart
%
%
\begin{equation}\label{eq:loss}
\widehat{\gamma}(\lambda)= - \frac{2}{n}\sum_{i=1}^n \mathsf{f}_\lambda
(X_i) + \| \mathsf{f}_\lambda\|^2.
\end{equation}
This motivates the use of $\widehat{\gamma} = \widehat{\gamma
}(\lambda)$ as the empirical criterion;
see, for instance, \cite{bm97,rud82,weg99}.

We define the
penalty
%
%
\begin{equation}\label{pen2}
\operatorname{pen}(\lambda)= 2 \sum
_{j=1}^M \omega_j |\lambda_j|
\end{equation}
with weights $\omega_j$ to be specified later, and we propose the
following data-driven
choice of $\lambda$:
%
%
\begin{eqnarray}\label{original}
\widehat{
\lambda}& = & \mathop{\arg\min}_{\lambda\in\mathbb R^M} \{
\widehat{\gamma}(\lambda)+ \operatorname{pen}(\lambda) \}
\nonumber\\[-8pt]\\[-8pt]
&=&
\mathop{\arg\min}_{\lambda\in\mathbb R^M} \Biggl\{ - \frac{2}{n}\sum
_{i=1}^n \mathsf{f}_\lambda(X_i) + \| \mathsf{f}_\lambda\|^2 + 2 \sum
_{j=1}^M \omega_j
|\lambda_j| \Biggr\}. \nonumber
\end{eqnarray}
Our estimator of
density $f$ that we will further call the \textit{SPADES estimator} is
defined by
\[
{f}^{\spadesuit}(x)=\mathsf{f}_{\widehat\lambda}(x)\qquad \forall x\in
\mathbb R^d.
\]
It is easy to see that, for an orthonormal system $\{f_j\}$, the
SPADES estimator coincides with the soft thresholding estimator
whose components are of the form $\widehat{
\lambda}_j=(1-\omega_j/|\tilde\lambda_j|)_+ \tilde\lambda_j$ where
$\tilde\lambda_j=n^{-1}\sum_{i=1}^n f_j(X_i)$ and $x_+ = \max(0,x)$.
We see that in this case $\omega_j$ is the threshold for the $j$th
component of a preliminary estimator
$\tilde\lambda=(\tilde\lambda_1,\ldots,\tilde\lambda_M)$.

The SPADES estimate can be easily computed by convex programming
even if $M \gg n$. We present an algorithm in Section \ref{sec6} below. SPADES
retains the desirable theoretical properties
of other density estimators, the computation of which may become
problematic for $M \gg n$. We refer to \cite{dl00} for a thorough
overview on combinatorial methods in density estimation, to
\cite{v99} for density estimation using support vector machines and
to \cite{bm97} for density estimates using penalties proportional to
the dimension.

\section{Oracle inequalities for SPADES}\label{sec3}

\subsection{Preliminaries}

For any $\lambda\in\mathbb R^M$, let
\[
J(\lambda)= \bigl\{ j\in\{1,\ldots,M\}\dvtx\lambda_j\ne0 \bigr\}
\]
be the set of indices corresponding to nonzero components of
$\lambda$ and
\[
M(\lambda)= | J(\lambda) | =\sum_{j=1}^M I\{ \lambda_j\ne0\}
\]
its cardinality. Here $I\{ \cdot\}$ denotes the indicator function.
Furthermore, set
\[
\sigma_j^2=\operatorname{Var}(f_j(X_1)),\qquad
L_j= \|f_j\|_\infty,
\]
for $1\le j\le M$, where $\operatorname{Var}(\zeta)$ denotes
the variance of random variable $\zeta$ and $\| \cdot\|_\infty$ is
the $L_\infty(\mathbb R^d)$ norm.

We will prove sparsity oracle inequalities for the estimator
$\widehat\lambda=\widehat\lambda(\omega_1,\ldots,\omega_M)$,
provided the
weights $\omega_j$ are chosen large enough. We first consider a
simple choice:
%
%
\begin{equation}
\label{a} \omega_j = 4 L_j r(\delta/2),
\end{equation}
where
$0<\delta<1$ is a user-specified parameter and
%
%
\begin{equation}\label{r}
r(\delta) = r(M,n,\delta) = \sqrt{ \frac{ \log(M/\delta) }{n} }.
\end{equation}
The oracle inequalities that we prove below hold with a probability
of at least $1-\delta$ and are nonasymptotic: they are valid for all
integers $M$ and $n$. The first of these inequalities is established
under a coherence condition on the ``correlations''
\[
\rho_M(i,j)= \frac{\langle f_i,f_j\rangle }{\|f_i\| \|f_j\|},\qquad i,j=1,\ldots,M.
\]
For $\lambda\in\mathbb R^M$, we define a local coherence number (called
\textit{maximal local coherence}) by
\[
\rho(\lambda)= \max_{i\in J(\lambda)} \max_{j\ne i}
|\rho_M(i,j)|,
\]
and we also define
\[
F(\lambda)= \max_{j\in J(\lambda)} \frac{\omega_j }{r(\delta/2) \|
f_j\|} = \max_{j\in J(\lambda)} \frac{4L_j }{\| f_j\|}
\]
and
\[
G=\max_{1\le j\le M} \frac{ r(\delta/2) \| f_j\|}{ \omega_j}
=\max_{1\le j\le M} \frac{ \| f_j\|}{ 4L_j} .
\]

\subsection{Main results}

\begin{theorem}\label{th:mutcoh}
Assume that $L_j<\infty$ for $1\le j\le M$. Then with probability at
least $1-\delta$ for all $\lambda\in\mathbb R^M$ that satisfy
%
%
\begin{equation}\label{mutcoh}
16 G F(\lambda) \rho(\lambda)
M(\lambda) \le1
\end{equation}
and all $\alpha>1$, we have the
following oracle inequality:
\begin{eqnarray*}
&&
\| {f}^{\spadesuit}-f\|^2 +\frac{\alpha}{2(\alpha-1)} \sum_{j=1}^M
\omega_j |\widehat\lambda_j-\lambda_j| \\
&&\qquad\le\frac{\alpha
+1}{\alpha-1} \|
\mathsf{f}_\lambda-f\|^2 + \frac{ 8\alpha^2 }{\alpha-1} F^2(\lambda)
r^2(\delta/2) M(\lambda).
\end{eqnarray*}
\end{theorem}

Note that only a condition on the local coherence (\ref{mutcoh}) is
required to obtain the result of Theorem \ref{th:mutcoh}. However,
even this condition can be too strong, because the bound on
``correlations'' should be \textit{uniform} over $j\in J(\lambda), i\ne
j$; cf. the definition of $\rho(\lambda)$. For example, this excludes
the cases where the ``correlations'' can be relatively large for a
small number of pairs $(i,j)$ and almost zero for otherwise. To
account for this situation, we suggest below another version of
Theorem \ref{th:mutcoh}. Instead of maximal local coherence, we
introduce \textit{cumulative local coherence} defined by
\[
\rho_*(\lambda)= \sum_{i\in J(\lambda)} \sum_{j> i}
|\rho_M(i,j)|.
\]
\begin{theorem}\label{th:cumcoh}
Assume that $L_j<\infty$ for $1\le j\le M$. Then with probability at
least $1-\delta$ for all $\lambda\in\mathbb R^M$ that satisfy
%
%
\begin{equation}\label{cumcoh}
16 F(\lambda) G \rho_*(\lambda)
\sqrt{M(\lambda)} \le1
\end{equation}
and all $\alpha>1$, we have
the following oracle inequality:
\begin{eqnarray*}
&&
\| {f}^{\spadesuit}-f\|^2 +\frac{\alpha}{2(\alpha-1)} \sum_{j=1}^M
\omega_j |\widehat\lambda_j-\lambda_j| \\
&&\qquad\le\frac{\alpha
+1}{\alpha-1} \|
\mathsf{f}_\lambda-f\|^2 + \frac{8\alpha^2}{\alpha-1} F^2(\lambda)
r^2(\delta/2) M(\lambda).
\end{eqnarray*}
\end{theorem}

Theorem \ref{th:cumcoh} is useful when we deal with sparse Gram
matrices $\Psi_M= ( \langle f_i,\break f_j\rangle  )_{1\le i,j\le M}$ that have
only a small number $N$ of nonzero off-diagonal entries. This
number will be called a \textit{sparsity index} of matrix $\Psi_M$, and
is defined as
\[
N= \bigl| \bigl\{(i,j)\dvtx i,j\in\{1,\ldots,M\}, i>j \mbox{ and }
\psi_M(i,j)\ne0 \bigr\}\bigr|,
\]
where $\psi_M(i,j)$ is the $(i,j)$th entry of $\Psi_M$ and $|A|$
denotes the cardinality of a set $A$. Clearly,
$N<M(M+1)/2$. We therefore obtain the following immediate corollary of
Theorem \ref{th:cumcoh}.
\begin{corollary} \label{cor:1} Let $\Psi_M$ be a Gram matrix with
sparsity index $N$. Then the assertion of Theorem \ref{th:cumcoh}
holds if we replace there (\ref{cumcoh}) by the condition
%
%
\begin{equation}\label{cumcoh1}
16 F(\lambda) N
\sqrt{M(\lambda)} \le1.
\end{equation}
\end{corollary}

We finally give an oracle inequality, which is valid under the
assumption that the Gram matrix $\Psi_M$ is positive definite. It is
simpler to use than the above results when the dictionary is
orthonormal or forms a frame. Note that the coherence assumptions
considered above do not necessarily imply the positive definiteness
of~$\Psi_M$. Vice versa, the positive definiteness of $\Psi_M$ does
not imply these assumptions.
\begin{theorem}\label{thm:1a} Assume that $L_j<\infty$ for
$1\le j\le M$ and that the Gram matrix $\Psi_M$ is positive definite
with minimal eigenvalue larger than or equal to $\kappa_M>0$. Then,
with probability at least $1-\delta$, for all $\alpha>1$ and all
$\lambda\in\mathbb R^M$, we have
%
%
\begin{eqnarray}\label{th11}
&&
\| {f}^{\spadesuit} - f\|^2 + \frac{\alpha}{\alpha-1}
\sum_{j=1}^M \omega_j
|\widehat{\lambda}_j-\lambda_j|\nonumber\\[-8pt]\\[-8pt]
&&\qquad \le
\frac{\alpha+1}{\alpha-1} \| \mathsf{f}_\lambda- f\|^2 +
\biggl(\frac{8\alpha^2} {\alpha-1} \biggr)\frac{ G(\lambda)}{n
\kappa_M},\nonumber
\end{eqnarray}
where
\[
G(\lambda)\triangleq\sum_{j\in J(\lambda)} \omega_j^2=
\frac{16 \log(2M/\delta)}{n}\sum_{j\in J(\lambda)} L_j^2 .
\]
\end{theorem}

We can consider some other choices for $\omega_j$ without
affecting the previous results. For instance,
%
%
\begin{equation}
\label{bb} \omega_j = 2\sqrt{2} \sigma_j r (\delta/2) + \tfrac83
L_j r^2(\delta/2)
\end{equation}
or
%
%
\begin{equation}
\label{c} \omega_j = 2\sqrt{2} T_j r(\delta/2) + \tfrac83L_j
r^2(\delta/2)
\end{equation}
with
\[
T_j^2 = \frac{2}{n} \sum_{i=1}^n f_j^2(X_i)
+ 2 L_j^2 r^2(\delta/2)
\]
yield the same conclusions. These modifications of (\ref{a}) prove
useful, for example, for situations where $f_j$ are wavelet basis
functions; cf. Section \ref{sec5}. The choice (\ref{c}) of $\omega_j$ has an
advantage of being completely data-driven.
\begin{theorem}\label{thm:3}
Theorems \ref{th:mutcoh}--\ref{thm:1a} and Corollary \ref{cor:1}
hold with the choices (\ref{bb}) or (\ref{c}) for the weights
$\omega_j$ without changing the assertions. They also remain valid
if we replace these $\omega_j$ by any $\omega_j'$ such that
$\omega_j'>\omega_j$.
\end{theorem}

If $\omega_j$ is chosen as in (\ref{c}), our bounds on the risk of
SPADES estimator involve the random variables $(1/n)\sum_{i=1}^n
f_j^2(X_i)$. These can be replaced in the bounds by deterministic
values using the following lemma.
\begin{lemma}\label{oneside}
Assume that $ L_j<\infty$ for $j=1,\ldots,M$. Then
%
%
\begin{equation}
\mathbb{P}\Biggl(\frac{1}{n} \sum_{i=1}^n f_j^2(X_i) \le2\mathbb{E}f_j^2(X_1)
+\frac43 L_j^2r^2(\delta/2), \forall j=1,\ldots,M \Biggr)\ge
1-\delta/2.\hspace*{-32pt}
\end{equation}
\end{lemma}

From Theorem \ref{thm:3} and Lemma \ref{oneside} we find
that, for the choice of $\omega_j$ as in (\ref{c}), the oracle
inequalities of Theorems \ref{th:mutcoh}--\ref{thm:1a} and Corollary
\ref{cor:1} remain valid with probability at least $1-3\delta/2$ if
we replace the $\omega_j$ in these inequalities by the expressions
$2\sqrt{2} \tilde T_j r(\delta/2) + (8/3)L_j r^2(\delta/2)$ where
$\tilde T_j= (2\mathbb{E}f_j^2(X_1)
+(4/3)L_j^2r^2(\delta/2) )^{1/2}$.

\subsection{Proofs}

We first prove the following preliminary lemma.
Define the random variables
\[
V_j = \frac{1}{n} \sum_{i=1}^n \{ f_j(X_i) - \mathbb{E}f_j(X_i) \}
\]
and the
event
%
%
\begin{equation}\label{A}
A = \bigcap_{j=1}^M \{ 2 | V_j| \le\omega_j \}.
\end{equation}
\begin{lemma}\label{prelim}
Assume that $ L_j<\infty$ for $j=1,\ldots,M$. Then for all
$\lambda\in\mathbb R^M$ we have that, on the event $A$,
%
%
\begin{equation}\qquad
\| {f}^{\spadesuit} -f\|^2 + \sum_{j=1}^M \omega_j
|\widehat{\lambda}_j-\lambda_j| \le\| \mathsf{f}_\lambda- f\|^2 + 4
\sum_{j\in J(\lambda)} \omega_j | \widehat{\lambda}_j -\lambda_j|.
\end{equation}
\end{lemma}
\begin{pf}
By the
definition of $\widehat{\lambda}$,
\[
- \frac{2}{n}\sum_{i=1}^n \mathsf{f}_{\widehat\lambda}(X_i) + \| \mathsf
{f}_{\widehat\lambda} \|^2 +
2 \sum_{j=1}^M \omega_j | \widehat\lambda_j| \le-
\frac{2}{n}\sum_{i=1}^n \mathsf{f}_\lambda(X_i) + \| \mathsf{f}_\lambda
\|
^2 + 2 \sum_{j=1}^M \omega_j |\lambda_j|
\]
for all $\lambda\in\mathbb R^M$. We rewrite this inequality as
\begin{eqnarray*}
\| {f}^{\spadesuit} -f \|^2 &\le& \| \mathsf{f}_\lambda-f \|^2 - 2\langle f,
{f}^{\spadesuit} - \mathsf{f}_\lambda\rangle  + \frac{2}{n} \sum_{i=1}^n
({f}^{\spadesuit} - \mathsf{f}_\lambda)(X_i)\\
&&{} + 2 \sum_{j=1}^M \omega_j
|\lambda_j| - 2 \sum_{j=1}^M \omega_j | \widehat\lambda_j| \\
&=& \| \mathsf{f}_\lambda-f \|^2 + 2\sum_{j=1}^M \Biggl( \frac{1}{n}
\sum_{i=1}^n f_j(X_i) - \mathbb{E}f_j(X_i) \Biggr)(\widehat\lambda_j -
\lambda_j) \\
&&{} + 2\sum_{j=1}^M \omega_j |\lambda_j| - 2
\sum_{j=1}^M \omega_j |\widehat\lambda_j|.
\end{eqnarray*}
Then, on the event $A$,
\[
\| {f}^{\spadesuit} - f\|^2 \le\| \mathsf{f}_\lambda-f \|^2 +
\sum_{j=1}^M \omega_j |\widehat\lambda_j - \lambda_j| + 2
\sum_{j=1}^M \omega_j |\lambda_j| - 2\sum_{j=1}^M \omega_j |
\widehat
\lambda_j|.
\]
Add $\sum_j \omega_j |\widehat{\lambda}_j
-\lambda_j|$ to both sides of the inequality to obtain
\begin{eqnarray*}
&& \| {f}^{\spadesuit} - f\|^2 + \sum_{j=1}^M\omega_j
|\widehat{\lambda}_j -\lambda_j| \\
&&\qquad \le\| \mathsf{f}_\lambda-f \|^2 +
2\sum_{j=1}^M \omega_j |\widehat\lambda_j - \lambda_j| + 2
\sum_{j=1}^M \omega_j |\lambda_j| -
2\sum_{j=1}^M \omega_j | \widehat\lambda_j|\\
&&\qquad
\le\| \mathsf{f}_\lambda- f\| ^2 +2
\sum_{j\in J(\lambda) } \omega_j |\widehat\lambda_j - \lambda_j|
+2 \sum_{j=1}^M \omega_j |\lambda_j| -2 \sum_{j\in J(\lambda)}
\omega_j
|\widehat{\lambda}_j|\\
&&\qquad
\le\| \mathsf{f}_\lambda- f\| ^2 +4
\sum_{j\in J(\lambda) } \omega_j |\widehat\lambda_j - \lambda_j|,
\end{eqnarray*}
where we used that $\lambda_j=0$ for $j\notin J(\lambda)$ and the
triangle inequality.
\end{pf}

For the choice (\ref{a}) for $\omega_j$, we find by Hoeffding's
inequality for sums of independent random variables $\zeta_{ij}=
f_j(X_i) - \mathbb{E}f_j(X_i)$ with $|\zeta_{ij}|\le2L_j$ that
\[
\mathbb{P}(A) \le\sum_{j=1}^M
\mathbb{P}\{ 2| V_j| > \omega_j \} \le
2 \sum_{j=1}^M \exp\biggl( - \frac{2n\omega_j^2/4}{8 L_j^2} \biggr) =\delta.
\]
\begin{pf*}{Proof of Theorem \ref{th:mutcoh}}
In view of Lemma \ref{prelim}, we need to bound $\sum_{j\in
J(\lambda)} \omega_j \times\break|\widehat\lambda_j -\lambda_j|$. Set
\[
u_j = \widehat\lambda_j -\lambda_j,\qquad U(\lambda) = \sum_{j\in
J(\lambda)} |u_j| \| f_j\|,\qquad U = \sum_{j=1}^M |u_j| \| f_j\|
r=r(\delta/2).
\]
Then, by the definition of $F(\lambda)$,
\[
\sum_{j\in J(\lambda)} \omega_j |\widehat\lambda_j -\lambda_j|
\le
r F(\lambda) U(\lambda).
\]
Since
\[
\mathop{\sum\sum}_{i,j\notin J(\lambda)}
\langle f_i,f_j\rangle  u_i u_j\ge0,
\]
we obtain
%
%
\begin{eqnarray}\label{Aa}\qquad
\sum_{j\in J(\lambda)} u_j^2 \|f_j\|^2 &=& \|{f}^{\spadesuit}-\mathsf
{f}_\lambda\|^2 -\mathop{\sum\sum}_{i,j\notin J(\lambda)} u_i u_j
\langle f_i,f_j\rangle
\nonumber\\
&&{}-
2\sum_{i\notin J(\lambda)} \sum_{j\in J(\lambda)} u_i u_j
\langle f_i,f_j\rangle -\mathop{\sum\sum}_{i,j\in J(\lambda), i\ne j} u_i u_j
\langle f_i,f_j\rangle \nonumber\\
&\le& \| {{f}^{\spadesuit}-\mathsf{f}_\lambda\|^2 + 2\rho(\lambda)
\sum_{i\notin
J(\lambda)}}| u_i| \| f_i\|\sum_{j\in J(\lambda)}| u_j |\| f_j\|
\\
&&{}+ \rho(\lambda)\mathop{\sum\sum}_{i,j\in J(\lambda)} | u_i| | u_j|
\|f_i\| \| f_j\|\nonumber\\
&=&\| {f}^{\spadesuit}-\mathsf{f}_\lambda\|^2 + 2\rho(\lambda) U
(\lambda) U - \rho(\lambda) U^2(\lambda).\nonumber
\end{eqnarray}
The left-hand side can be bounded by $\sum_{j\in J(\lambda)} u_j^2
\|f_j\|^2 \ge U^2(\lambda)/M(\lambda)$ using the Cauchy--Schwarz
inequality, and we obtain that
\[
U^2(\lambda) \le\|{f}^{\spadesuit}-\mathsf{f}_\lambda\|^2 M(\lambda) +
2\rho(\lambda) M(\lambda) U(\lambda) U,
\]
which immediately implies
%
%
\begin{equation}\label{B}
U(\lambda) \le2\rho(\lambda) M(\lambda) U + \sqrt{M(\lambda)}
\|{f}^{\spadesuit} -\mathsf{f}_\lambda\|.
\end{equation}
Hence, by Lemma \ref{prelim}, we have, with probability at least
$1-\delta$,
\begin{eqnarray*}
&& \| {f}^{\spadesuit} -f\|^2 + \sum_{j=1}^M \omega_j
|\widehat{\lambda}_j-\lambda_j| \\
&&\qquad\le \| \mathsf{f}_\lambda- f\|^2 + 4
\sum_{j\in J(\lambda)} \omega_j | \widehat{\lambda}_j -\lambda
_j|\\
&&\qquad\le
\| \mathsf{f}_\lambda- f\|^2 +4r F(\lambda) U(\lambda)\\
&&\qquad\le \| \mathsf{f}_\lambda- f\|^2 + 4r F(\lambda) \bigl\{2\rho(\lambda)
M(\lambda) U + \sqrt{M(\lambda)} \| {f}^{\spadesuit}-\mathsf{f}_\lambda
\| \bigr\}\\
&&\qquad\le \| \mathsf{f}_\lambda- f\|^2 + 8 F(\lambda) \rho(\lambda)
M(\lambda) G \sum_{j=1}^M \omega_j |\widehat{\lambda}_j -\lambda
_j| \\
&&\qquad\quad{} +
4r F(\lambda)\sqrt{M(\lambda)} \| {f}^{\spadesuit}-\mathsf{f}_\lambda\| .
\end{eqnarray*}
For all $\lambda\in\mathbb R^M$ that satisfy relation (\ref
{mutcoh}), we
find that, with probability exceeding $1-\delta$,
\begin{eqnarray*}
&&
\| {f}^{\spadesuit} -f\|^2 + \frac12 \sum_{j=1}^M \omega_j
|\widehat{\lambda}_j-\lambda_j|
\\
&&\qquad\le \|\mathsf{f}_\lambda- f\|^2 +
4 r
F(\lambda) \sqrt{M(\lambda)} \| {f}^{\spadesuit} -\mathsf{f}_\lambda
\| \\
&&\qquad\le \| \mathsf{f}_\lambda- f\|^2 + 2 \bigl\{ 2 r F(\lambda)
\sqrt{M(\lambda)} \bigr\} \| {f}^{\spadesuit}- f \| \\
&&\qquad\quad{}+ 2
\bigl\{ 2 r F(\lambda) \sqrt{M(\lambda)} \bigr\}
\| \mathsf{f}_\lambda-f \|.
\end{eqnarray*}
After applying the inequality $2xy\le x^2/\alpha+ \alpha y^2$
($x,y\in\mathbb R, \alpha>0$) for each of the last two summands, we
easily find
the claim.
\end{pf*}
\begin{pf*}{Proof of Theorem \ref{th:cumcoh}}
The proof is similar to that of Theorem
\ref{th:mutcoh}. With
\[
U_*(\lambda)=\sqrt{ \sum_{j\in J(\lambda)} u_j^2 \|f_j\|^2},
\]
we obtain now the following analogue of (\ref{Aa}):
\begin{eqnarray*}
U_*^2(\lambda) &\le& {\| {f}^{\spadesuit}-\mathsf{f}_\lambda\|^2 +
2\rho_*(\lambda) \max_{i\in J(\lambda), j>i}}| u_i| \| f_i\|
| u_j |\| f_j\|\\
&\le& \| {f}^{\spadesuit}-\mathsf{f}_\lambda\|^2 +
2\rho_*(\lambda)U_*(\lambda) \sum_{j=1}^M | u_j| \| f_j\|
\\
&=&\| {f}^{\spadesuit}-\mathsf{f}_\lambda\|^2 + 2\rho_*(\lambda) U_*
(\lambda) U.
\end{eqnarray*}
Hence, as in the proof of Theorem \ref{th:mutcoh}, we have
\[
U_*(\lambda) \le2\rho_*(\lambda) U + \|{f}^{\spadesuit} -\mathsf
{f}_\lambda\|,
\]
and using the inequality $U_*(\lambda)\ge
U(\lambda)/\sqrt{M(\lambda)}$, we find
%
%
\begin{equation}\label{B1}
U(\lambda) \le2\rho_*(\lambda)\sqrt{M(\lambda)} U +
\sqrt{M(\lambda)} \|{f}^{\spadesuit} -\mathsf{f}_\lambda\|.
\end{equation}
Note that (\ref{B1}) differs from (\ref{B}) only in the fact that
the factor $2\rho(\lambda)M(\lambda)$ on the right-hand side is now
replaced by $2\rho_*(\lambda)\sqrt{M(\lambda)}$. Up to this
modification, the rest of the proof is identical to that of Theorem
\ref{th:mutcoh}.
\end{pf*}
\begin{pf*}{Proof of Theorem \ref{thm:1a}}
By the assumption on $\Psi_M$, we have
\[
\| \mathsf{f}_\lambda\|^2 = \mathop{\sum\sum}_{1\le i,j\le M} \lambda
_i\lambda_j
\int_{\mathbb R^d} f_i(x) f_j(x) \,dx \ge\kappa_M \sum_{j\in
J(\lambda)} \lambda_j^2.
\]
By the Cauchy--Schwarz inequality, we find
\begin{eqnarray*}
&&
4 \sum_{j\in J(\lambda)} \omega_j | \widehat{\lambda}_j - \lambda_j|
\\
&&\qquad\le 4 \sqrt{\sum_{j\in J(\lambda)} \omega_j^2} \sqrt{\sum
_{j\in
J(\lambda)}
|\widehat{\lambda}_j -\lambda_j|^2} \\
&&\qquad\le 4 \biggl(\frac{\sum_{j\in J(\lambda)}
\omega_j^2}{n\kappa_M} \biggr)^{1/2} \| {f}^{\spadesuit}- \mathsf{f}_\lambda\|.
\end{eqnarray*}
Combination with Lemma \ref{prelim} yields that, with probability at
least $1-\delta$,
%
%
\begin{eqnarray}\label{simple}
&&\| {f}^{\spadesuit} -f\|^2 + \sum_{j=1}^M \omega_j
|\widehat{\lambda}_j-\lambda_j| \nonumber\\
&&\qquad\le \| \mathsf{f}_\lambda- f\|^2 + 4
\biggl(\frac{\sum_{j\in J(\lambda)}
\omega_j^2}{n\kappa_M} \biggr)^{1/2} \|
{f}^{\spadesuit}- \mathsf{f}_\lambda\| \\
&&\qquad\le \| \mathsf{f}_\lambda- f\|^2 + b ( \| {f}^{\spadesuit}-f\| +
\|\mathsf{f}_\lambda-f\| ) ,\nonumber
\end{eqnarray}
where $b=4 \sqrt{\sum_{j\in J(\lambda)} \omega_j^2}/ \sqrt{
n\kappa_M}$. Applying the inequality $2xy\le x^2/\alpha+ \alpha y^2$
($x,y\in\mathbb R, \alpha>0$) for each of the last two summands in
(\ref{simple}), we get the result.
\end{pf*}
\begin{pf*}{Proof of Theorem \ref{thm:3}}
Write $\bar\omega_j = 2\sqrt{2} \sigma_j r (\delta/2) + (8/3) L_j
r^2(\delta/2)$ for the choice of $\omega_j$ in (\ref{bb}). Using
Bernstein's exponential inequality for sums of independent random
variables $\zeta_{ij}= f_j(X_i) - \mathbb{E}f_j(X_i)$ with
$|\zeta_{ij}|\le2L_j$, we obtain that
%
%
\begin{eqnarray}\label{bern}
\mathbb{P}(A^c) &=& \mathbb{P}\Biggl( \bigcup_{j=1}^M \{ 2|V_j|>
\bar\omega_j\} \Biggr)
\nonumber\\
&\le&\sum_{j=1}^M \mathbb{P}\{ 2|V_j|>
\bar\omega_j\}\nonumber\\[-8pt]\\[-8pt]
&\le& \sum_{j=1}^M
\exp\biggl( -\frac{n \bar\omega_j^2/4}{2\operatorname{Var}(f_j(X_1)) + 2
L_j \bar\omega_j /3} \biggr)\nonumber\\
&\le& M\exp\bigl( -nr^2(\delta/2)\bigr) = \delta/2.\nonumber
\end{eqnarray}
Let now $\omega_j$ be defined by (\ref{c}). Then, using
(\ref{bern}), we can write
%
%
\begin{eqnarray}\label{bern1}
\mathbb{P}(A^c) &=& \mathbb{P}\Biggl( \bigcup_{j=1}^M \{ 2|V_j|> \omega
_j\} \Biggr)
\nonumber\\
&\le& \sum_{j=1}^M \mathbb{P}\{ 2|V_j|> \bar\omega_j\} +
\sum_{j=1}^M \mathbb{P}\{ \bar\omega_j > \omega_j\}
\\
&\le& \delta/2 +\sum_{j=1}^M \mathbb{P}\{ \bar\omega_j > \omega
_j\} .\nonumber
\end{eqnarray}
Define
\[
t_j =2 \frac{ \mathbb{E}f_j^4(X_1)}{ \mathbb{E}f_j^2(X_1)} \frac
{\log
(2M/\delta)}{n}
\]
and note that
\[
\frac{2}{n} \sum_{i=1}^n f_j^2(X_i) + t_j \le T_j^2.
\]
Then
\begin{eqnarray*}
\mathbb{P}\{\bar\omega_j > \omega_j\} &=&
\mathbb{P}\{ \operatorname{Var}(f_j(X_1)) > T_j^2 \}\\
&\le& \mathbb{P}\Biggl\{ \mathbb{E}f_j^2(X_1) > \frac{2}{n}\sum_{i=1}^n
f_j^2(X_i) + t_j \Biggr\} \\
&\le& \exp\biggl(- \frac{ n \{ \mathbb{E}f_j^2(X_1) + t_j \}^2}{8 \mathbb
{E}f_j^4(X_1)} \biggr)\qquad
\mbox{using Proposition 2.6 in \cite{weg03}}\\
&\le& \exp\biggl( -\frac{nt_j \mathbb{E}f_j^2(X_1) }{2 \mathbb{E}f
_j^4(X_1) } \biggr)
\qquad \mbox{since } (x+y)^2 \ge4xy,
\end{eqnarray*}
which is less than $\delta/(2M)$.
Plugging this in (\ref{bern1}) concludes the proof.
\end{pf*}
\begin{pf*}{Proof of Lemma \ref{oneside}} Using
Bernstein's\vspace*{1pt} exponential inequality for sums of independent random
variables $f_j^2(X_i) - \mathbb{E}f_j^2(X_i)$ and the fact that
$\mathbb{E}
f_j^4(X_1) \le L_j^2 \mathbb{E}f_j^2(X_1)$, we find
\begin{eqnarray*}
&&\mathbb{P}\Biggl(\frac{1}{n} \sum_{i=1}^n f_j^2(X_i) \ge2\mathbb{E}f_j^2(X_1)
+\frac43 L_j^2r^2(\delta/2) \Biggr)
\\
&&\qquad= \mathbb{P}\Biggl(\frac{1}{n} \sum_{i=1}^n f_j^2(X_i) - \mathbb{E}f_j^2(X_1)
\ge\mathbb{E}f_j^2(X_1) +\frac43 L_j^2r^2(\delta/2) \Biggr)
\\
&&\qquad\le \exp\biggl( -\frac{n (\mathbb{E}f_j^2(X_1) +4/3
L_j^2r^2(\delta/2))^2}{2 \mathbb{E}f_j^4(X_1) + 4/3 L_j^2\{
\mathbb{E}
f_j^2(X_1) +4/3 L_j^2r^2(\delta/2)\} } \biggr)
\\
&&\qquad\le \exp\bigl( -nr^2(\delta/2)\bigr) = \frac{\delta}{2M} ,
\end{eqnarray*}
which implies the lemma.
\end{pf*}

\section{Sparse estimation in mixture models}\label{sec4}

In this section we assume that the true density $f$ can be
represented as a finite mixture
\[
f(x)=\sum_{j \in I^*}{\bar\lambda}_j p_j(x),
\]
where $I^* \subseteq\{1, \ldots, M\}$ is unknown, $p_j$ are known
probability densities and ${\bar\lambda}_j> 0$ for all $j\in I^*$.
We focus in this section on model selection, that is, on the correct
identification of the set $I^*$. It will be convenient for us to
normalize the densities $p_j$ by their $L_2$ norms and to write the
model in the form
\[
f(x)=\sum_{j \in I^*}\lambda_j^* f_j(x),
\]
where $I^* \subseteq\{1, \ldots, M\}$ is unknown, $f_j=p_j/\|p_j\|$
are known functions and $\lambda_j^*>0 $ for all $j\in I^*$. We set
$\lambda^*=(\lambda_1^*,\ldots,\lambda_M^*)$, where $\lambda_j^*=0,
j\notin I^*$.

For clarity of exposition, we consider a simplified version of the
general setup introduced above. We compute the estimates of
$\lambda^*$ via (\ref{original}), with weights defined by [cf.
(\ref{a})]:
\[
\omega_j = 4Lr\qquad \mbox{for all } j,
\]
where $r>0$ is a constant that we specify below, and for clarity of
exposition we replaced all $L_j =\| f_j\|_\infty$ by an upper bound
$L$ on $\max_{1 \leq j \leq M}L_j$. Recall that, by construction,
$\|f_j\| = 1$ for all $j$. Under these assumptions condition
(\ref{mutcoh}) takes the form
%
%
\begin{equation}\label{sim} \rho(\lambda) \leq
\frac{1}{16M(\lambda)} .
\end{equation}
We state (\ref{sim}) for the true vector $\lambda^*$ in the
following form:
\renewcommand{\theCondition}{(\Alph{Condition})}
\begin{Condition}\label{condA}
\[
\rho^* \leq\frac{1}{16k^*},
\]
where $k^* = |I^*| = M(\lambda^*)$ and $\rho^*=\rho(\lambda^*)$.
\end{Condition}

Similar conditions are quite standard in the literature on sparse
regression estimation and compressed sensing; cf., for example,
\cite{det04,zy,btw05,btw06b,brt,bun08}. The difference is that those
papers use the empirical version of the correlation $\rho^*$ and the
numerical constant in the inequality is, in general, different from
$1/16$. Note that Condition \ref{condA} is quite intuitive. Indeed, the
sparsity index $k^*$ can be viewed as the effective dimension of the
problem. When $k^*$ increases the problem becomes harder, so that we
need stronger conditions (smaller correlations $\rho^*$) in order to
obtain our results. The interesting case that we have in mind is
when the effective dimension $k^*$ is small, that is, the model is
sparse.

The results of Section \ref{sec3} are valid for any $r$ larger or
equal to $r(\delta/2)= \{ \log(2M/\delta) / {n} \}^{1/2}$. They give
bounds on
the predictive performance of\break SPADES. As noted in, for example,
\cite{bun08}, for $\ell_1$-penalized model selection in regression,
the tuning sequence $\omega_j$ required for correct selection is
typically larger than the one that yields good prediction. We show
below that the same is true for selecting the components of a
mixture of densities. Specifically, in this section we will take the
value
%
%
\begin{equation}\label{newr} r = r\bigl(M, n, \delta/(2M)\bigr) =
\sqrt{ \frac{ \log(2M^2/\delta) }{n} }.
\end{equation}
We will use the following corollary of Theorem
\ref{th:mutcoh}, obtained for $\alpha= \sqrt{2}$.
\begin{corollary}\label{cormixt} Assume that
Condition \ref{condA} holds. Then with probability at least
$1-\delta/M$, we have
%
%
\begin{equation}\label{buna}
\sum_{j=1}^{M}|\widehat{\lambda}_j - \lambda_j^*| \leq
\frac{4\sqrt{2}}{L}k^*\sqrt{ \frac{\log(2M^2/\delta)}{n} } .
\end{equation}
\end{corollary}

Inequality
(\ref{buna}) guarantees that the estimate $\widehat{\lambda}$ is
close to the true $\lambda^*$ in $\ell_1$ norm, if the number of
mixture components $k^*$ is substantially smaller than $\sqrt{n}$.
We regard this as an intermediate step for the next result that
deals with the identification of $I^*$.

\subsection{Correct identification of the mixture components}\label{sec41}

We now show that $I^*$ can be identified with probability close to 1
by our procedure. Let $\widehat{I}=J(\widehat{\lambda})$ be the set of
indices of the nonzero components of $\widehat{\lambda}$ given by
(\ref{original}). In what follows we investigate when $P(\widehat{I} =
I^*) \geq1-\varepsilon$ for a given $0 < \varepsilon< 1$. Our
results are nonasymptotic, they hold for any fixed $M$ and $n$.

We need two conditions to ensure that correct recovery of $I^*$ is
possible. The first one is the identifiability of the model, as
quantified by Condition \ref{condA} above. The second condition
requires that the weights of the mixture are above the noise level,
quantified by $r$. We state it as follows:
\begin{Condition}\label{condB}
\[
\min_{j \in I^*} |\lambda_j^*| > 4\bigl(\sqrt{2}+1\bigr)r L,
\]
where $L=\max(1/\sqrt{3}, \max_{ 1\leq j \leq M}L_j )$ and
$r$ is given in (\ref{newr}).
\end{Condition}
\begin{theorem}\label{whole} Let $ 0 < \delta< 1/2$ be a given number.
Assume that Conditions \ref{condA} and \ref{condB} hold. Then $\mathbb{P}(\widehat{I} =I^*)
\geq1 - 2\delta(1 + 1/M)$.
\end{theorem}
\begin{remark*} Since all $\lambda_j^*$ are nonnegative, it seems
reasonable to restrict the minimization in (\ref{original}) to
$\lambda$ with nonnegative components. Inspection of the proofs
shows that all the results of this section remain valid for such a
modified estimator. However, in practice, the nonnegativity issue is
not so important. Indeed, the estimators of the weights are
quite close to the true values and turn out to be positive for
positive $\lambda_j^*$. For example, this was the case in our
simulations discussed in Section \ref{sec6} below. On the other hand, adding
the nonnegativity constraint in (\ref{original}) introduces some
extra burden on the numerical algorithm. More generally, it is
trivial to note that the results of this and previous sections
extend verbatim to the setting where $\lambda\in\Lambda$ with
$\Lambda$ being any subset of ${\mathbb R}^M$. Then the minimization in
(\ref{original}) should be performed on $\Lambda$, in the theorems
of Section \ref{sec3} we should replace $\lambda\in\mathbb R^M$ by $\lambda
\in
\Lambda$ and in this section $\lambda^*$ should be supposed to
belong to $\Lambda$.
\end{remark*}
\begin{pf*}{Proof of Theorem \ref{whole}}
We begin by noticing
that
\[
\mathbb{P}(\widehat{I} \neq I^*) \leq\mathbb{P}(I^* \not\subseteq
\widehat{I}) +
\mathbb{P}(\widehat{I} \not\subseteq I^*),
\]
and we control each of the probabilities on the right-hand side separately.

\textit{Control of $\mathbb{P}(I^* \not\subseteq\widehat{I})$}. By
the definitions of the sets $\widehat{I}$ and $I^*$, we have
\begin{eqnarray*}
\mathbb{P}(I^* \not\subseteq\widehat{I})
&\leq& \mathbb{P}(\widehat\lambda_k = 0 \mbox{ for some } k\in I^*)
\\
& \leq& k^*\max_{k \in I^*} \mathbb{P}(\widehat\lambda_k = 0 ).
\end{eqnarray*}
We control the last probability by using the characterization
(\ref{cond}) of $\widehat{\lambda}$ given in Lemma \ref{sol} of the
\hyperref[app]{Appendix}. We also recall that $\mathbb{E}f_k(X_1) = \sum_{j \in I^*}
\lambda_j^*\langle f_k, f_j \rangle= \sum_{j=1}^M
\lambda_j^*\langle f_k, f_j \rangle$, since we assumed that the
density of $X_1$ is the mixture $f^* = \sum_{j \in I^*}
\lambda_j^*f_j$. We therefore obtain, for $k \in I^*$,
%
%
\begin{eqnarray}
\mathbb{P} ( \widehat\lambda_k = 0 )
&=& \mathbb{P} \Biggl( \Biggl|\frac{1}{n}\sum_{i=1}^{n}f_k(X_i) - \sum
_{j=1}^{M}\widehat{\lambda}_j\langle f_j, f_k\rangle\Biggr|
\leq4rL; \widehat\lambda_k = 0 \Biggr) \nonumber\\
& = & \mathbb{P} \Biggl( \Biggl|\frac{1}{n}\sum_{i=1}^{n}f_k(X_i) -
\mathbb{E}f_k(X_1)\nonumber\\
&&\hspace*{18pt}{} - \sum_{j=1}^{M}(\widehat{\lambda}_j -
\lambda_j^*)\langle f_j, f_k\rangle\Biggr| \leq4rL; \widehat\lambda_k = 0 \Biggr)
\nonumber\\
& \le& \mathbb{P} \Biggl( \Biggl| \lambda_k^*\|f_k\|^2 +
\frac{1}{n}\sum_{i=1}^{n}f_k(X_i) - \mathbb{E}f_k (X_1)\nonumber\\
&&\hspace*{68.1pt}{}- \sum_{j
\neq k}(\widehat{\lambda}_j - \lambda_j^*)\langle f_j, f_k\rangle\Biggr|
\leq4rL \Biggr) \nonumber\\
\label{zeroo}
& \leq& \mathbb{P} \Biggl( \Biggl| \frac{1}{n}\sum_{i=1}^{n}f_k(X_i) - \mathbb
{E}f_k(X_1) \Biggr| \geq\frac{|\lambda_k^*|\|f_k\|^2}{2} - 2rL \Biggr) \\
\label{zero}
&&{} + \mathbb{P} \biggl( \biggl| \sum_{j\ne k}(\widehat{\lambda}_j -
\lambda_j^*)\langle f_j, f_k\rangle\biggr|\geq
\frac{|\lambda_k^*|\|f_k\|^2}{2} - 2rL \biggr).
\end{eqnarray}
To bound (\ref{zeroo}), we use Hoeffding's inequality, as in the
course of Lemma \ref{prelim}. We first recall that $ \|f_k\| = 1$
for all $k$ and that, by Condition \ref{condB}, $ \min_{k \in
I^*}|\lambda_k^*| \geq4(\sqrt{2}+1)Lr$, with $ r = r(\delta/(2M))=
\{{\log(2M^2/\delta)}/{n}\}^{1/2}$. Therefore,
%
%
\begin{eqnarray} \label{unu} && \mathbb{P}
\Biggl( \Biggl| \frac{1}{n}\sum_{i=1}^{n}f_k(X_i) - \mathbb{E}f_k(X_1) \Biggr|
\geq\frac{ |\lambda_k^*|}{2} - 2rL \Biggr) \nonumber\\[-8pt]\\[-8pt]
&&\qquad \leq\mathbb{P}
\Biggl( \Biggl| \frac{1}{n}\sum_{i=1}^{n}f_k(X_i) - \mathbb{E}f_k(X_1) \Biggr| \geq
2\sqrt{2}rL \Biggr)
\leq\frac{\delta}{M^2}. \nonumber
\end{eqnarray}
To bound
(\ref{zero}), notice that, by Conditions \ref{condA} and \ref{condB},
\begin{eqnarray*}
&& \mathbb{P} \biggl( \biggl| \sum_{j \ne k}(\widehat{\lambda}_j - \lambda
_j^*)\langle f_j, f_k\rangle\biggr|\geq\frac{|\lambda_k^*|}{2} - 2rL \biggr)
\\
&&\qquad \leq
\mathbb{P} \Biggl( \sum_{j =1}^{M}|\widehat{\lambda}_j - \lambda_j^*|
\geq32\sqrt{2}rLk^* \Biggr) \leq\mathbb{P} \Biggl( \sum_{j =1}^{M}|\widehat{\lambda}_j - \lambda
_j^*| \geq\frac{4\sqrt{2}rk^*}{L} \Biggr) \leq\frac{\delta}{M},
\end{eqnarray*}
where the penultimate inequality holds since, by
definition, $L^2 \ge
1/3$ and the last inequality holds by Corollary \ref{cormixt}.

Combining the above results, we obtain
\[
\mathbb{P}(I^* \not\subseteq\widehat{I}) \leq k^*\frac{\delta}{M^2}
+ k^*\frac{\delta}{M} \leq\frac{\delta}{M} + \delta.
\]

\textit{Control of $\mathbb{P}(\widehat{I} \not\subseteq I^*)$.} Let
%
%
\begin{equation}\label{reduce}
h(\mu) = -\frac2n
\sum_{i=1}^n \sum_{j \in I^*}\mu_jf_j(X_{i}) + \biggl\| \sum_{j \in
I^*}\mu_jf_j\biggr\|^2 +
8rL \sum_{j \in I^*} |\mu_j|.
\end{equation}
Let
%
%
\begin{eqnarray}\label{(2.9)}
\tilde{\mu} = \mathop{\arg\min}_{\mu\in\mathbb R^{k^{*}}} h(\mu).
\end{eqnarray}
Consider the random event
%
%
\begin{equation}\label{b}
\mathcal{B} = \bigcap_{k \notin I^*} \Biggl\{
\Biggl|-\frac{1}{n}\sum_{i=1}^{n}f_k(X_{i}) + \sum_{ j \in I^*}
\tilde{\mu}_j\langle f_j, f_k \rangle\Biggr| \le4L r \Biggr\}.
\end{equation}
Let $\bar{\mu} \in\mathbb R^M$ be the vector that has the components of
$\tilde{\mu}$ given by (\ref{(2.9)}) in positions corresponding to the
index set $I^*$ and zero components elsewhere. By the first part of
Lemma \ref{sol} in\vspace*{1pt} the
\hyperref[app]{Appendix},\vspace*{1pt} we have that $\bar{\mu} \in
\mathbb R^M$ is a solution of (\ref{original}) on the event
$\mathcal{B}$. Recall that $\widehat{\lambda}$ is also a solution of
(\ref{original}). By the definition of the set $\widehat{I}$, we have
that $\widehat{\lambda}_k \neq0$ for $k \in\widehat{I}$. By
construction, $\tilde{\mu}_k \neq0$ for some subset $S \subseteq I^*$.
By the second part of Lemma \ref{sol} in the \hyperref[app]{Appendix},
any two solutions have nonzero elements in the same positions.
Therefore, $\widehat{I} = S \subseteq I^*$ on $\mathcal{B}$. Thus,
%
%
\begin{eqnarray}\label{eq:main}
\mathbb{P}(\widehat{I} \not\subseteq I^*) &\leq& \mathbb{P}( \mathcal
{B}^c) \nonumber\\
&\leq&\sum_{k \notin I^*} \mathbb{P}
\Biggl\{ \Biggl|-\frac{1}{n}\sum_{i=1}^{n}f_k(X_{i}) + \sum_{ j \in
I^*} \tilde{\mu}_j\langle f_j, f_k \rangle\Biggr| \ge
4rL \Biggr\} \nonumber\\[-8pt]\\[-8pt]
&\leq&\sum_{k \notin I^*}\mathbb{P} \Biggl( \Biggl|\frac{1}{n}\sum
_{i=1}^{n}f_{k}(X_i) - Ef_k(X_1) \Biggr| \geq2\sqrt{2}rL \Biggr) \nonumber\\
&&{} + \sum_{k \notin I^*} \mathbb{P} \biggl(\sum_{j \in
I^*}|\tilde{\mu}_j - \lambda^*_j| |\langle f_j, f_k \rangle
| \geq\bigl(4-2\sqrt{2}\bigr)rL \biggr). \nonumber
\end{eqnarray}
Reasoning as in (\ref{unu}) above, we find
\[
\sum_{k \notin I^*}\mathbb{P} \Biggl( \Biggl|\frac{1}{n}\sum_{i=1}^{n}f_{k}(X_i)
- Ef_k(X_1) \Biggr| \geq2\sqrt{2}rL \Biggr) \leq\frac{\delta}{M}.
\]

To bound the last sum in (\ref{eq:main}), we first notice that
Theorem \ref{th:mutcoh} [if we replace there $r(\delta/2)$ by the
larger value $r(\delta/(2M))$; cf. Theorem
\ref{thm:3}] applies to $\tilde{\mu}$ given by~(\ref{(2.9)}). In
particular,
\[
\mathbb{P} \biggl(\sum_{j \in I^*}|\tilde{\mu}_j - \lambda^*_j| \geq
\frac{4\sqrt{2}}{L}k^*r \biggr) \leq\frac{\delta}{M}.
\]
Therefore, by Condition \ref{condA}, we have
\begin{eqnarray*}
&&\sum_{k \notin I^*} \mathbb{P} \biggl(\sum_{j \in I^*}
|\tilde{\mu}_j - \lambda^*_j| |\langle f_j, f_k \rangle|\geq
\bigl(4-2\sqrt{2}\bigr)rL \biggr) \\
&&\qquad \leq\sum_{k \notin I^*} \mathbb{P} \biggl(\sum_{j \in I^*}
|\tilde{\mu}_j - \lambda^*_j| \geq32\bigl(4-2\sqrt{2}\bigr)k^*rL \biggr)
\\
&&\qquad \leq\sum_{k \notin I^*} \mathbb{P} \biggl(\sum_{j \in
I^*}|\tilde{\mu}_j - \lambda^*_j| \geq
\frac{4\sqrt{2}}{L}k^*r \biggr) \leq\delta,
\end{eqnarray*}
which holds since $L^2 \geq1/3$. Collecting all the
bounds above, we obtain
\[
P(\widehat{I} \neq I^*) \leq2\delta+ \frac{2\delta}{M},
\]
which concludes the proof.
\end{pf*}

\subsection{Example: Identifying true components in mixtures of
Gaussian densities}\label{sec42}

Consider an ensemble of $M$ Gaussian densities $p_j$'s in $\mathbb R^d$ with
means $\mathbb\mu_j$ and covariance matrices $\tau_j\mathbb I_d$, where
$\mathbb I_d$ is the unit $d \times d$ matrix. In what follows we show
that Condition \ref{condA} holds if the means of the Gaussian
densities are well separated and we make this precise below.
Therefore, in this case, Theorem \ref{whole} guarantees that if the
weights of the mixture are above the threshold given in Condition \ref{condB},
we can recover the true mixture components with high probability via
our procedure. The densities are
\[
p_j(x) =
\frac{1}{(2\pi\tau_j^2)^{d/2} }\exp\biggl(-\frac{\|x -
\mu_j\|_2^2}{2\tau_j^2} \biggr),
\]
where \mbox{$\| \cdot\|_2$}
denotes the Euclidean norm. Consequently, $f_j= p_j/\| p_j\|$ with
$\| p_j\|= (4\pi\tau_j^2)^{-d/4}$. Recall that
Condition \ref{condA} requires
\[
16\rho^* = 16 \max_{i \in I^*, j \neq i} |
\langle f_i,f_j\rangle | \leq1/k^*.
\]
Let $\tau_{\max} = \max_{ 1 \leq j \leq
M}\tau_j$ and $D_{\min}^{2} = \min_ {k \neq j} \|\mu_k -
\mu_j\|_2^2$. Via simple algebra, we obtain
\[
\rho^* \leq\exp\biggl( - \frac{D_{\min}^2}{4\tau_{\max}^2} \biggr).
\]
Therefore, Condition \ref{condA} holds if
%
%
\begin{equation}\label{condmu}
D_{\min}^2 \geq4\tau_{\max}^2\log(16k^*).
\end{equation}
Using this and Theorem \ref{whole}, we see that SPADES
identifies the true components in a mixture of Gaussian densities if
the square Euclidean distance between any two means is large enough
as compared to the largest variance of the components in the
mixture.

Note that Condition \ref{condB} on the size of the mixture weights
involves the constant~$L$, which in this example can be taken as
\[
L= \max
\biggl( \frac{\sqrt{3} }{3},
\max_{ 1 \leq j \leq M}\|f_j\|_{\infty} \biggr)
=
\max\biggl(
\frac{\sqrt{3} }{3} , ( \pi\tau_{\min}^2 )^{-d/4} \biggr),
\]
where $\tau_{\min} = \min_{ 1 \leq j \leq M}\tau_j$.
\begin{remark*}
Often both the location and scale parameters are unknown. In this
situation, as suggested by the Associate Editor, the SPADES procedure
can be
applied to a family of densities with both scale and location
parameters chosen from an appropriate grid. By Theorem \ref
{th:mutcoh}, the resulting estimate will be a good approximation of
the unknown target density. An immediate modification of Theorem \ref
{whole}, as in \cite{bun08b}, further guarantees that SPADES
identifies correctly the important components
of this approximation.
\end{remark*}

\section{SPADES for adaptive nonparametric density estimation}\label{sec5}

We assume in this section that the density $f$ is defined on a
bounded interval of $\mathbb R$ that we take without loss of generality to
be the interval $[0,1]$. Consider a countable system of
functions\vspace*{1pt}
$\{\psi_{lk}, l\ge-1, k\in V(l)\}$ in $L_2$, where the set of
indices $V(l)$ satisfies $|V(-1)|\le C$, $2^l\le|V(l)|\le C2^l,
l\ge0$, for some constant $C$, and where the functions $\psi_{lk}$
satisfy
%
%
\begin{equation}\label{snorm}
\|\psi_{lk}\|\le C_1,\qquad \|\psi_{lk}\|_\infty\le C_1 2^{l/2},\qquad
\biggl\|\sum_{k\in V(l)}\psi_{lk}^2 \biggr\|_\infty\le C_1 2^{l}
\end{equation}
for all $l\ge-1$ and for some $C_1<\infty$. Examples of such systems
$\{\psi_{lk}\}$ are
given, for instance, by compactly supported wavelet bases; see,
for example, \cite{hkpt98}. In this case $\psi_{lk}(x) =
2^{l/2}\psi(2^lx-k)$ for some compactly supported function $\psi$.
We assume that
$\{\psi_{lk}\}$ is a frame, that is, there exist positive constants
$c_1$ and $c_2$ depending only on $\{\psi_{lk}\}$ such that, for any
two sequences of coefficients $\beta_{lk}$, $\beta_{lk}'$,
%
%
\begin{eqnarray}\label{frame}
c_1\sum_{l=-1}^\infty\sum_{k\in V(l)}(\beta_{lk}-\beta_{lk}')^2
&\le&
\Biggl\|\sum_{l=-1}^\infty\sum_{k\in
V(l)}(\beta_{lk}-\beta_{lk}')\psi_{lk} \Biggr\|^2 \nonumber\\[-8pt]\\[-8pt]
&\le& c_2
\sum_{l=-1}^\infty\sum_{k\in V(l)}(\beta_{lk}-\beta_{lk}')^2.\nonumber
\end{eqnarray}
If $\{\psi_{lk}\}$ is an orthonormal wavelet basis, this condition
is satisfied with \mbox{$c_1=c_2=1$}.

Now, choose $\{f_1,\ldots,f_M\} = \{\psi_{lk}, -1\le l \le
l_{\max}, k\in V(l)\}$, where $l_{\max}$ is such that
$2^{l_{\max}}\asymp n/(\log n)$. Then also $M\asymp n/(\log n)$. The
coefficients $\lambda_j$ are now indexed by $j=(l,k)$, and we set by
definition $\lambda_{(l,k)}=0$ for $(l,k)\notin\{-1\le l \le
l_{\max}, k\in V(l)\}$. Assume that there exist coefficients
$\beta_{lk}^*$ such that
\[
f= \sum_{l=-1}^\infty\sum_{k\in V(l)}\beta_{lk}^*\psi_{lk},
\]
where the series converges in $L_2$. Then Theorem \ref{thm:1a}
easily implies the following result.
\begin{theorem}\label{th:adnonp}
Let $f_1,\ldots,f_M$ be as defined above with $M\asymp n/(\log n)$,
and let $\omega_j$ be given by (\ref{c}) for $\delta=n^{-2}$. Then
for all $n\ge1$, $\lambda\in\mathbb R^M$ we have, with probability
at least
$1-n^{-2}$,
%
%
\begin{eqnarray}\label{th:adnonp1}\quad
\| {f}^{\spadesuit} - f\|^2 & \le&
K \Biggl( \sum_{l=-1}^\infty\sum_{k\in
V(l)}\bigl(\lambda_{(l,k)}-\beta_{lk}^*\bigr)^2 \nonumber\\[-8pt]\\[-8pt]
&&\hspace*{14.2pt}{}+\sum_{(l,k)\in J(\lambda)}
\Biggl[\frac1{n}\sum_{i=1}^n \psi_{lk}^2(X_i)\frac{\log n}{n} +
2^l \biggl(\frac{\log n}{n} \biggr)^2 \Biggr] \Biggr),\nonumber
\end{eqnarray}
where $K$ is a constant independent of $f$.
\end{theorem}

This is a general oracle inequality that allows one to show that the
estimator ${f}^{\spadesuit}$ attains minimax rates of convergence,
up to a logarithmic factor simultaneously on various functional
classes. We will explain this in detail for the case where $f$ belongs
to a
class of functions ${\mathcal F}$ satisfying the following
assumption for some $s>0$:
\begin{Condition}\label{condC}
For any $f\in{\mathcal F}$ and any
$l'\ge0$ there exists a sequence of coefficients
$\lambda=\{\lambda_{(l,k)}, -1\le l \le l', k\in V(l) \}$ such that
%
%
\begin{equation}\label{appro}
\sum_{l=-1}^\infty\sum_{k\in
V(l)}\bigl(\lambda_{(l,k)}-\beta_{lk}^*\bigr)^2\le C_2 2^{-2l's}
\end{equation}
for a constant $C_2$ independent of $f$.
\end{Condition}

It is well known that Condition \ref{condC} holds for various functional
classes ${\mathcal F}$, such as H\"{o}lder, Sobolev, Besov classes, if
$\{\psi_{lk}\}$ is an appropriately chosen wavelet basis; see, for example,
\cite{hkpt98} and the references cited therein. In this case $s$ is
the smoothness parameter of the class. Moreover, the basis
$\{\psi_{lk}\}$ can be chosen so that Condition \ref{condC} is satisfied
with $C_2$ independent of $s$ for all $s\le s_{\max}$, where
$s_{\max}$ is a given positive number. This allows for adaptation in
$s$.

Under Condition \ref{condC}, we obtain from (\ref{th:adnonp1}) that, with
probability at least $1-n^{-2}$,
%
%
\begin{eqnarray}\label{th:adnonp2}\qquad
\| {f}^{\spadesuit} - f\|^2 &\le&\min_{l'\le l_{\max}}
K \Biggl( C_2 2^{-2l's} +\sum_{(l,k)\dvtx l\le l'}
\Biggl[\frac1{n}\sum_{i=1}^n \psi_{lk}^2(X_i)\frac{\log n}{n}\nonumber\\[-8pt]\\[-8pt]
&&\hspace*{157.9pt}{} +
2^l \biggl(\frac{\log n}{n} \biggr)^2 \Biggr] \Biggr).\nonumber
\end{eqnarray}
From (\ref{th:adnonp2}) and the last inequality in (\ref{snorm}) we
find for some constant $K'$, with probability at least $1-n^{-2}$,
%
%
\begin{eqnarray}\label{th:adnonp3}
\| {f}^{\spadesuit} - f\|^2 & \le&\min_{l'\le l_{\max}}
K' \biggl( 2^{-2l's} + 2^{l'} \biggl(\frac{\log n}{n} \biggr) +
2^{2l'} \biggl(\frac{\log n}{n} \biggr)^2 \biggr) \nonumber\\[-8pt]\\[-8pt]
&=&
O \biggl( \biggl(\frac{\log n}{n} \biggr)^{-2s/(2s+1)}
\biggr),\nonumber
\end{eqnarray}
where the last expression is obtained by choosing $l'$ such that
$2^{l'}\asymp(n/\break\log n)^{1/(2s+1)}$. It follows from
(\ref{th:adnonp3}) that ${f}^{\spadesuit}$ converges with the
optimal rate (up to a logarithmic factor) simultaneously on all the
functional classes satisfying Condition \ref{condC}. Note that the
definition of the functional class is not used in the construction
of the estimator ${f}^{\spadesuit}$, so this estimator is optimal
adaptive in the rate of convergence (up to a logarithmic factor) on
this scale of functional classes for $s\le s_{\max}$.
Results of such type, and even more pointed (without extra
logarithmic factors in the rate and sometimes with exact asymptotic
minimax constants), are known for various other adaptive density
estimators; see, for instance, \cite{g92,bm97,hkpt98,kpt96,r06,rt04}
and the references cited therein. These papers consider classes of
densities that are uniformly bounded by a fixed constant; see the
recent discussion in~\cite{bir08}. This prohibits, for example, free
scale transformations of densities within a class. Inequality
(\ref{th:adnonp3}) does not have this drawback. It allows to get the
rates of convergence for classes of unbounded densities $f$ as well.

Another example is given by the classes of sparse densities defined
as follows:
\[
{\mathcal L}_0(m)= \bigl\{f\dvtx[0,1]\to\mathbb R\dvtx\mbox{$f$ is a
probability density and } |\{j\dvtx\langle f,f_j\rangle  \ne0\} |\le
m \bigr\},
\]
where $m\le M$ is an unknown integer. If $f_1,\ldots,f_M$ is a
wavelet system as defined above and $J^*=\{j=(l,k)\dvtx\langle f,f_j\rangle  \ne
0\}$, then under the conditions of Theorem \ref{th:adnonp} for any
$f\in{\mathcal L}_0(m)$ we have, with probability at least
$1-n^{-2}$,
%
%
\begin{equation}\label{th:adnonp5}
\| {f}^{\spadesuit} - f\|^2 \le
K \Biggl( \sum_{(l,k)\in J^*} \Biggl[\frac1{n}\sum_{i=1}^n
\psi_{lk}^2(X_i)\frac{\log n}{n} + 2^l \biggl(\frac{\log
n}{n} \biggr)^2 \Biggr] \Biggr) .
\end{equation}
From (\ref{th:adnonp5}), using Lemma \ref{oneside} and the first two
inequalities in (\ref{snorm}), we obtain the following result.
\begin{corollary} \label{cor:twee}
Let the assumptions of Theorem \ref{th:adnonp} hold. Then, for every
$L<\infty$ and $n\ge1$,
%
%
\begin{eqnarray} \label{eq:cor:twee}
\sup_{f\in{\mathcal L}_0(m)\cap\{f\dvtx\|f\|_\infty\le
L\}}\mathbb{P}\biggl\{\|{f}^{\spadesuit} -f\|^2 \ge b \biggl(\frac{m\log
n}{n} \biggr) \biggr\} \le(3/2)n^{-2}\nonumber\\[-8pt]\\[-8pt]
\eqntext{\forall m\le M,}
\end{eqnarray}
where
$b>0$ is a constant depending only on $L$.
\end{corollary}

Corollary \ref{cor:twee} can be viewed as an analogue for density
estimation of the
adaptive minimax results for ${\mathcal L}_0$ classes obtained in
the Gaussian sequence model \cite{abdj,g02} and in the random design
regression model \cite{btw06b}.

\section{Numerical experiments}\label{sec6}

In this section we describe the algorithm used for the minimization
problem (\ref{original}) and we assess the performance of our
procedure via a simulation study.

\subsection{A coordinate descent algorithm}\label{sec61}
Since the criterion given in (\ref{original}) is convex, but not
differentiable, we adopt an optimization by coordinate descent
instead of a gradient-based approach (gradient descent, conjugate
gradient, etc.) in the spirit of \cite{fhht07,fht08}. Coordinate
descent is an iterative greedy optimization technique that starts at
an initial location $\lambda\in\mathbb R^M$ and at each step chooses one
coordinate $\lambda_j \in\mathbb R$ of $\lambda$ at random or in order
and finds the optimum in that direction, keeping the other variables
$\lambda_{-j}$ fixed at their current values. For convex functions,
it usually converges to the global optimum; see \cite{fhht07}. The method
is based on the obvious observation that for functions of the type
\[
H(\lambda) = g(\lambda)+\omega|\lambda|_1,
\]
where $g$ is a generic convex and differentiable function,
$\omega>0$ is a given parameter, and $|\lambda|_1$ denotes the
$\ell_1$ norm, the optimum in a direction $\lambda_j \in\mathbb R$
is to
the left, right or at $\lambda_j=0$, depending on the signs of the
left and right partial derivatives of $H$ at zero. Specifically, let
$g_j$ denote the partial derivative of $g$ with respect to
$\lambda_j$, and denote by
$\lambda^0_{-j}$ the vector $\lambda$ with the $j$th coordinate set
to 0. Then, the minimum in direction $j$ of $H(\lambda)$ is at
$\lambda^0_{-j}$ if and only if
$|g_j(\lambda^0_{-j}) |< \omega$. This observation makes the
coordinate descent become the iterative thresholding algorithm
described below.

\subsubsection*{Coordinate descent}

\mbox{}

Given $\omega$, initialize all $\lambda_j$, $1 \leq j \leq M$, for
example, with $1/M$.
\begin{enumerate}
\item Choose a direction $j\in\{1,\ldots,M\}$ and set $\lambda
^{\mathrm{old}}=\lambda$.
\item If $|g_j(\lambda^0_{-j})| < \omega$, then set $\lambda
=\lambda^0_{-j}$, otherwise obtain $\lambda$ by
line minimization in direction $j$.
\item If $|\lambda^{\mathrm{old}}-\lambda|>\epsilon$, go to 1, where
$\epsilon>0$ is a given precision level.
\end{enumerate}
For line minimization, we used the procedure \texttt{linmin} from
Numerical Recipes \cite{nr07}, page 508.

\subsection{Estimation of mixture weights using the generalized
bisection method and a penalized cross-validated loss function}
\label{sec:gbm}

We apply the coordinate descent algorithm described above to
optimize the function $H(\lambda)$ given by (\ref{original}), where
the tuning parameters $\omega_j$ are all set to be equal to the same
quantity $\omega$. The theoretical choice of this quantity described
in detail in the previous sections may be too conservative in
practice. In this section we propose a data driven method for
choosing the tuning parameter $\omega$, following the procedure
first introduced in \cite{bb09}, which we briefly describe here for
completeness.

The procedure chooses adaptively the tuning
parameter from a list of candidate values, and it has two
distinctive features: the list of candidates is not given by a fine
grid of values and the adaptive choice is not given by
cross-validation, but by a dimension stabilized cross-validated criterion.
We begin by describing the principle underlying our construction of
the set of candidate values which, by avoiding a grid search, provides
significant computational savings. We use a generalization of the bisection
method to find, for each $0\leq k\leq M$, a preliminary tuning
parameter $\omega= w_k$ that gives a solution $\widehat{\lambda}^k$
with exactly $k$ nonzero elements. Formally, denote by $\widehat
n(\omega)$ the number of nonzero elements in the $\lambda$ obtained
by minimizing (\ref{original}) with $\omega_j\equiv\omega$
for a given value of the tuning parameter $\omega$. The generalized
bisection method will find a sequence of values of the tuning
parameter, $w_0,\ldots,w_M$, such that $\widehat n(w_k)=k$, for each $0 \leq
k\leq M$. It proceeds as follows, using a queue consisting of pairs
$(w_i, w_j)$ such that $\widehat n(w_i)<\widehat n(w_j) - 1$.

\subsubsection*{The general bisection method (GBM) for all $k$}

\mbox{}

Initialize all $w_i$ with $-1$.
\begin{enumerate}
\item Choose $w_0$ very large, such that $\widehat n(w_0) = 0$. Choose
$w_n = 0$,
hence, $\widehat n(w_n) = n$.
\item Initialize a queue $q$ with the pair $(w_0, w_n)$.
\item Pop the first pair $(a, b)$ from the queue.
\item Take $w = (a+b)/2$. Compute $k = \widehat n(w)$.
\item If $w_k = -1$, make $w_k = w$.
\item If $|\widehat n(a) -k|>1$ and $|a -w|>\alpha$, add $(a, w)$ to
the back of the queue.
\item If $|\widehat n(b )-k|>1$ and $|b-w|>\alpha$, add $(\omega,
b)$ to the back of the queue.
\item If the queue is not empty, go to 3.
\end{enumerate}

This algorithm generalizes the basic bisection method (BBM),
which is a well-established computationally efficient method for
finding a root $z \in\mathbb R$ of a function $h(z)$; see, for example,
\cite{bf01}. We experimentally observed (see also \cite{bb09} for a
detailed discussion) that using the GBM is about 50 times faster
than a grid search with the same accuracy.

Our procedure finds the final tuning parameter $\omega$ by combining
the GBM with the dimension stabilized $p$-fold cross-validation
procedure summarized below. Let $D$ denote the whole data set, and let
$D=D_1\cup\cdots\cup D_p$ be a partition of $D$ in $p$ disjoint
subsets. Let $D_{-j} = D \setminus D_j$. We will denote by $w_k^j$ a
candidate tuning parameter determined using the GBM on $D_{-j}$.
We denote by $I_k^j$ the set of indices corresponding to the nonzero
coefficients of the estimator of $\lambda$ given by (\ref{original}),
for tuning parameter $w_k^j$ on $D_{-j}$. We denote by $ \widehat{\lambda
}^{kj}$ the minimizers on
$D_{-j}$ of the unpenalized criterion $\widehat\gamma(\lambda)$,
with respect only to those $\lambda_l$ with $l \in I_k^j$. Let $L_k^j
=: \widehat{\gamma}(\widehat{\lambda}^{kj})$, computed on $D_j$. With
this notation, the procedure becomes the following:

\subsubsection*{Weight selection procedure}

\mbox{}

Given: a data set $D$ partitioned into $p$ disjoint subsets, $D=D_1\cup
\cdots\cup D_p$. Let $D_{-j} = D \setminus D_j$ for all $j$.
\begin{enumerate}
\item For each $ 1 \leq k \leq M$ and each fold $j$ of the
partition, $ 1 \leq j \leq p$:

Use the GBM to find $w_k^j$ and $I_k^j$ such that $\widehat n(w_k^j)
=|I_k^j|=k$ on $D_{-j}$.

Compute $L_k^j=: \widehat{\gamma}(\widehat{\lambda}^{kj})$, as defined
above, on $D_j$.
\item For each $1 \leq k \leq M$:

Compute $L_k =: \frac{1}{p}\sum_{j = 1}^{p}L_k^j$.
\item Obtain
\[
\widehat{k} = \mathop{\arg\min}_{k} \biggl(L_k + 0.5k \frac{\log n}{n}\biggr).
\]
\item With $\widehat{k}$ from Step 3, use the BBM on the whole
data set $D$ to find the tuning sequence $w_{\widehat{k}}$ and then
compute the final estimators using the coordinate descent algorithm and
tuning paramemter $\omega= w_{\widehat{k}}$.
\end{enumerate}

In all the the numerical experiments described below we
took the number of splits $p=10$.
\begin{remark*}
We recall that the theoretical results of Section \ref{sec41}
show that for correct identification of the mixture components one
needs to work with a value of the tuning sequence that is slightly
larger than the one needed for good approximations with mixtures of a
given density. A good practical approximation of the latter tuning
value is routinely obtained by cross-validation; this approximation is,
however, not appropriate if the goal is correct selection, when the
theoretical results indicate that a different value is needed. Our
modification of the cross-validated loss function via a BIC-type
penalty is motivated by the known properties of the BIC-type criteria
to yield consistent model selection in a large array of models; see,
for example, \cite{bu04} for results on regression models.
The numerical experiments presented below show that this is also the
case for our criterion in the context of selecting mixture components.
The theoretical investigation of this method is beyond the scope of
this paper and will be undertaken in future research.
\end{remark*}

\subsection{Numerical results}\label{sec63}

In this subsection we illustrate the performance of our procedure via a
simulation study.

\subsubsection{One-dimensional densities} We begin by investigating
the ability of SPADES, with its tuning parameter chosen as above,
to (i) approximate well, with respect to the $L_2$ norm, a true
mixture; (ii) to identify the true mixture components. We conducted
a simulation study where the true density is a mixture of
Gaussian densities with $k^*=2$ and, respectively, $k^* = 5$ true
mixture components. The mixture components are chosen at random
from a larger pool of $M$ Gaussians ${\mathcal N}(aj,1)$, $ 1 \leq j
\leq M$, where for $k^*=2$ we take $a=4$, and for $k^*=5$ we take
$a=5$. These choices for $a$ ensure that the identifiability
condition (\ref{condmu}) is satisfied. The true components
correspond to the first $k^*$ Gaussian densities from our list, and
their weights in the true mixture are all equal to $1/k^*$. The
maximum size $M$ of the candidate list we considered is $M= 200$,
for $k^* = 2$ and $M = 600$, for $k^* = 5$. All the results obtained
below are relative to $S = 100$ simulations. Each time, a sample of
size $n$ is obtained from the true mixture and is the input of the
procedure described in Section \ref{sec:gbm}.

%
\begin{figure}[b]

\includegraphics{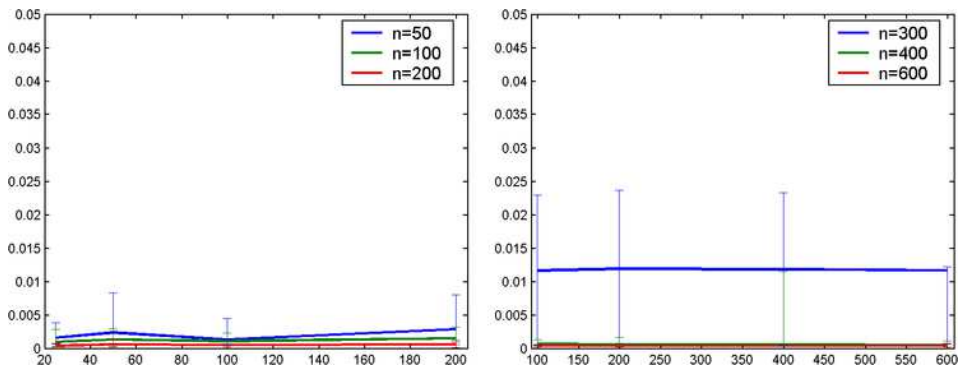}

\caption{Median $L_2$ error $\Vert f^*- f^{\spadesuit}\Vert^2$
for $|I^*|=2$, respectively, $|I^*|=5$. The error bars are the 25 and
75 percentiles.} \label{fig:err_m}
\end{figure}

We begin by evaluating the accuracy with respect to the $L_2$ norm
of the estimates of $f^*$. We investigate the sensitivity of our
estimates relative to an increase in the dictionary size and $k^*$.
In Figure \ref{fig:err_m}, we plot the median over 100 simulations
of $\|f^*- f^{\spadesuit}\|^2$ versus the size $M$ of the
dictionary, when the true mixture cardinality is $k^* = 2$ (left
panel) and $k^* = 5$ (right panel). For $k^* = 2$ we considered
three instances of sample sizes $n = 50, 100, 200$ and we varied $M$
up to 200. For $k^* = 5$ we considered three larger instances of sample
sizes $n = 300, 400, 600$ and we varied $M$ up to 600. These
experiments provide strong support for our theoretical results: the
increase in $M$ does not significantly affect the quality of
estimation, and an increase in $k^*$ does. For larger values of
$k^*$ we need larger
sample sizes to obtain good estimation accuracy.

%
\begin{figure}

\includegraphics{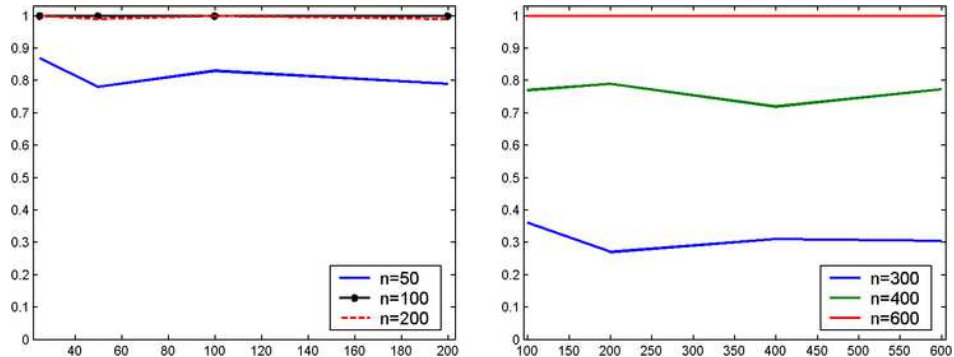}

\caption{Percentage of times $I^*=\widehat I$ obtained from
100 runs, for $|I^*|=2$, respectively, $|I^*|=5$.}
\label{fig:eye_m}
\end{figure}

We next investigated the ability of the SPADES to find the exact
mixture components. Figure \ref{fig:eye_m} shows a plot of the
percentage of times the exact mixture components were found versus
$M$. We considered the same combinations $(n, M)$ as in Figure
\ref{fig:err_m}. Again, observe that the performance does not
seriously degrade with the dictionary size $M$, and is almost
unaffected by its increase once a threshold sample size is being
used. However, notice that on the difference from the results
presented in Figure \ref{fig:err_m}, correct identification is poor
below the threshold sample size, which is larger for larger $k^*$.
This is in accordance with our theoretical results: recall
Condition \ref{condB} of Section \ref{sec41} on the minimum size of the mixture
weights. Indeed, we designed our simulations so that the weights are
relatively small for $k^* = 5$, they are all equal to $1/k^*=0.2$,
and a larger sample size is needed for their correct identification.

%
\begin{figure}[b]

\includegraphics{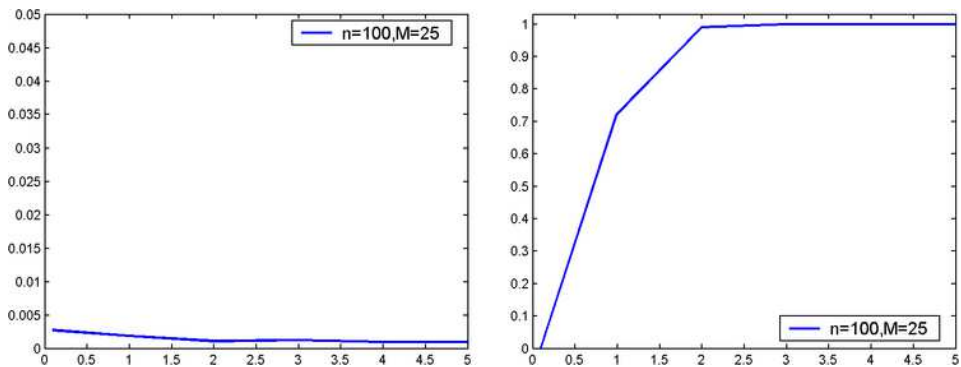}

\caption{Dependence on the distance $D_{\min} =
\min_{k \neq j} |\mu_k-\mu_j| $ of the $L_2$ error $\|f^*-
f^{\spadesuit} \|^2$ and the percentage of times $I^*=\widehat I$. In
this example, $n=100,M=25, |I^*|=2$.} \label{fig:err_d}
\end{figure}

Finally,\vspace*{2pt} we evaluated in Figure \ref{fig:err_d} the dependence of
the error and hit rate (i.e., the percentage of times $I^*=\widehat I$)
on the smallest distance $D_{\min} = \min_{k \neq j}
|\mu_k-\mu_j|$ between the means of the densities in the dictionary.
The results presented in Figures \ref{fig:err_m} and
\ref{fig:eye_m} above were obtained for the value $D_{\min} =
4$, which satisfies the theoretical requirement for correct mixture
identification. On the other hand, $D_{\min}$~can be smaller
for good $L_2$ mixture approximation. It is interesting to see what
happens when $D_{\min}$ decreases, so that the mixture
elements become very close to one another. In Figure \ref{fig:err_d}
we present the simulations for $k^* = 2$, $M = 25$ and $n = 100$,
which is sufficient to illustrate this point. We see that, although
the $L_2$ error increases slightly when $D_{\min} $ decreases,
the deterioration is not crucial. However, as our theoretical
results suggest, the percentage of times we can correctly identify
the mixture decreases to zero when the dictionary functions are very
close to each other.

\subsubsection{Two-dimensional densities}
\begin{figure}[b]

\includegraphics{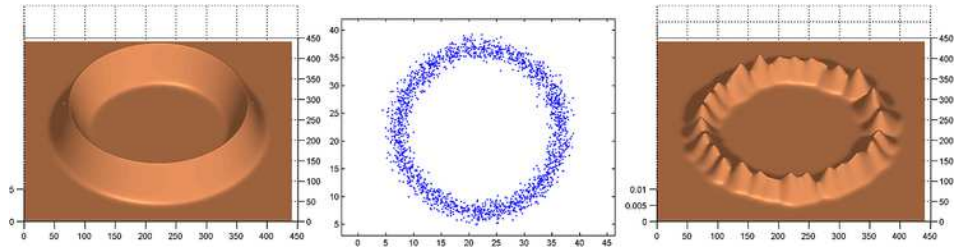}

\caption{A thick circle density, a sample of size 2000
from this density and approximations using a mixture of 80 isotropic
Gaussians.} \label{fig:plot_rec}
\end{figure}

In a second set of experiments our aim was to approximate a
two-dimensional probability density on a thick circle (cf. the left
panel of Figure \ref{fig:plot_rec}) with a mixture of isotropic
Gaussians. A sample of size 2000 from the circle density is shown in
the middle panel of Figure \ref{fig:plot_rec}. We use a set of
isotropic Gaussian candidates with covariance $\Sigma= \mathbb I_2$
centered at some of the 2000 locations, such that the Euclidean
distance between the means of any two such Gaussians is at least 1.
We select from these candidate mixture densities in a greedy
iterative manner, each time choosing one of the 2000 locations that
is at distance at least 1 from each of those already chosen. As a
result, we obtain a dictionary of $M=248$ candidate densities.

The circle density cannot be exactly represented as a finite mixture
of Gaussian components. This is a standard instance of many
practical applications in Computer Vision, as the statistics of
natural images are highly kurtotic and cannot be exactly
approximated by isotropic Gaussians. However, in many practical
applications a good approximation of an object that reflects its
general shape is sufficient and constitutes a first crucial step in
any analysis. We show below that SPADES offers such an
approximation.

Depending on the application, different trade-offs between the number
of mixture components (which relates to the computational demand of
the mixture model) and accuracy might be appropriate. For example,
in real-time applications a small number of mixture elements would
be required to fit into the computational constraints of the system,
as long as there is no significant loss in accuracy.

\begin{figure}

\includegraphics{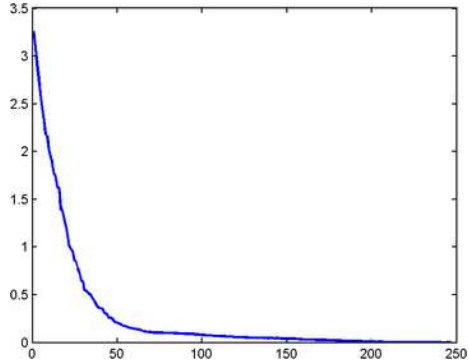}

\caption{Plot of $\widehat{\gamma}(\widehat
\lambda^k)-\gamma_0$ as a function of the mixture components $k$.}
\label{fig:loss}
\end{figure}

For the example presented below we used the GBM to determine
the mixture weights $\widehat\lambda^k$, for mixtures with $ k = 1, 2,
\ldots, 248$ components. Let
$\gamma_0 = \min_{k}\widehat{\gamma}(\widehat\lambda^k)$, where
we recall that the loss function $\widehat{\gamma}$ is given by
(\ref{eq:loss}) above. We used the quantity $\widehat{\gamma
}(\widehat\lambda^k)-\gamma_0$ to
measure the accuracy of the mixture approximation. In Figure
\ref{fig:loss}
%
%
we plotted $\widehat{\gamma}(\widehat\lambda
^k)-\gamma_0$ as
a function of $k$ and used this plot to determine the desired
trade-off between accuracy and mixture complexity. Based on this
plot, we selected the number of mixture components to be 80; indeed,
including more components does not yield any significant
improvement. The obtained mixture is displayed in the right panel of
Figure \ref{fig:plot_rec}.
%
%
We see that it successfully approximates
the circle density with a relatively small number of components.

\begin{appendix}\label{app}
\section*{Appendix}

\begin{lemma}\label{sol}
\textup{(I)} Let $\tilde{\mu}$ be given by (\ref{(2.9)}). Then
$\bar\mu=(\tilde{\mu}, 0) \in\mathbb R^M$ is a minimizer in
$\lambda\in\mathbb R^M$ of
\[
g(\lambda) = - \frac{2}{n}\sum_{i=1}^n \mathsf{f}_\lambda(X_i) + \|
\mathsf{f}_\lambda\|^2
+ 8Lr\sum_{k=1}^{M}|\lambda_k|
\]
on the random event $\mathcal{B}$ defined in (\ref{b}).

\textup{(II)} Any two minimizers of $g(\lambda)$ have nonzero
components in the same positions.
\end{lemma}
\begin{pf}
(I). Since $g$ is convex, by standard results in
convex analysis, $\bar{\lambda} \in\mathbb R^M$ is a minimizer of
$g$ if
and only if $0 \in D_{\bar{\lambda}}$ where $D_{\lambda}$ is the
subdifferential of $g(\lambda)$:
\begin{eqnarray*}
D_{\lambda} &=& \Biggl\{ w \in\mathbb R^M\dvtx w_k =
-\frac{2}{n}\sum_{i=1}^{n}f_k(X_{i}) + 2\sum_{ j =1}^M
\lambda_j\langle f_j, f_k \rangle+ 8rv_k,\\
&&\hspace*{154.2pt} v_k\in V_k(\lambda_k), 1
\leq k \leq M \Biggr\},
\end{eqnarray*}
where
\[
V_k(\lambda_k) = \cases{
\{L\}, &\quad if $\lambda_k > 0$, \cr
\{-L\}, &\quad if $\lambda_k < 0$, \cr
[-L, L], &\quad if $\lambda_k = 0$.}
\]
Therefore,
$\bar{\lambda}$ minimizes $g(\cdot)$
if and only if, for all $ 1 \leq k \leq M$,
%
%
\begin{eqnarray}
\label{cond}
\frac{1}{n}\sum_{i=1}^{n}f_k(X_{i}) - \sum_{ j =1}^M
\bar{\lambda}_j\langle f_j, f_k \rangle&= & 4Lr
\operatorname{sign}(\bar{\lambda}_k)\qquad \mbox{if }
\bar{\lambda}_k \neq0, \\
\label{condd}
\Biggl|\frac{1}{n}\sum_{i=1}^{n}f_k(X_{i}) - \sum_{ j =1}^M\bar{
\lambda}_j\langle f_j, f_k \rangle\Biggr| &\leq& 4Lr\qquad
\mbox{if } \bar{\lambda}_k = 0.
\end{eqnarray}
We now show that $\bar{\mu} = (\tilde{\mu}, 0) \in\mathbb R^M$ with
$\tilde{\mu}$ given in (\ref{(2.9)}) satisfies
(\ref{cond}) and (\ref{condd}) on the event ${\mathcal{B}}$ and
therefore is a minimizer of $g(\lambda)$ on this event. Indeed,
since $\tilde{\mu}$ is a minimizer of the convex function $h(\mu)$
given in (\ref{reduce}), the same convex analysis argument as above
implies that
\begin{eqnarray*}
\frac{1}{n}\sum_{i=1}^{n}f_k(X_{i}) - \sum_{ j \in I^*}
\tilde{\mu}_j\langle f_j, f_k \rangle
&= & 4Lr \operatorname{sign}(\tilde{\mu}_k)\qquad \mbox{if } \tilde{\mu}_k
\neq0, k \in I^*, \\
\Biggl| \frac{1}{n}\sum_{i=1}^{n}f_k(X_{i}) - \sum_{ j \in I^*}
\tilde{\mu}_j\langle f_j, f_k \rangle\Biggr| & \leq& 4Lr
\qquad\mbox{if } \tilde{\mu}_k = 0 , k \in I^*.
\end{eqnarray*}
Note that on the event ${\mathcal{B}}$ we also have
\begin{eqnarray}
\Biggl| \frac{1}{n}\sum_{i=1}^{n}f_k(X_{i}) - \sum_{ j \in I^*}
\tilde{\mu}_j\langle f_j, f_k \rangle\Biggr| \leq4Lr\nonumber\\
\eqntext{\mbox{if } k \notin I^*\ (\mbox{for which } \bar{\mu}_k = 0,
\mbox{by construction}).}
\end{eqnarray}
Here $\bar{\mu}_k$ denotes the $k$th coordinate of $\bar{\mu}$. The
above three displays and the fact that $\bar{\mu}_k=\tilde{\mu}_k, k
\in I^*$, show that $\bar{\mu}$ satisfies conditions
(\ref{cond}) and (\ref{condd}) and is therefore a minimizer of
$g(\lambda)$ on the event $\mathcal{B}$.

(II). We now prove the second assertion of the
lemma. In view of (\ref{cond}), the index set $S$ of the nonzero
components of any minimizer $\bar{\lambda}$ of $g(\lambda)$
satisfies
\[
S = \Biggl\{ k \in\{1, \ldots, M\}\dvtx
\Biggl|\frac{1}{n}\sum_{i=1}^{n}f_k(X_{i}) - \sum_{ j =1}^M
\bar{\lambda}_j\langle f_j, f_k \rangle\Biggr| = 4rL \Biggr\}.
\]
Therefore, if for any two minimizers $\bar{\lambda}^{(1)}$ and
$\bar{\lambda}^{(2)}$ of $g(\lambda)$ we have
%
%
\begin{equation}\label{same1}
\sum_{ j =1}^M \bigl(\bar{\lambda}_j^{(1)} - \bar{\lambda}_j^{(2)}\bigr)
\langle f_j, f_k \rangle= 0\qquad \mbox{for all } k,
\end{equation}
then $S$ is the same for all minimizers of $g(\lambda)$.

Thus, it remains to show (\ref{same1}). We use simple properties of
convex functions. First, we recall that the set of minima of a
convex function is convex. Then, if $\bar{\lambda}^{(1)}$ and
$\bar{\lambda}^{(2)}$ are two distinct points of minima, so is
$\rho\bar{\lambda}^{(1)} + (1 - \rho)\bar{\lambda}^{(2)}$, for
any $
0 < \rho< 1$. Rewrite this convex combination as
$\bar{\lambda}^{(2)} + \rho\eta$, where $\eta= \bar{\lambda}^{(1)}
- \bar{\lambda}^{(2)}$. Recall that the minimum value of any convex
function is unique. Therefore, for any $ 0 <\rho< 1$, the value
of $g(\lambda)$ at $\lambda=\bar{\lambda}^2 + \rho\eta$ is equal to
some constant $C$:
\begin{eqnarray*}
F(\rho)& \triangleq& -\frac{2}{n}
\sum_{i=1}^{n}\sum_{j=1}^{M} \bigl(\bar{\lambda}^{(2)}_j +
\rho\eta_j \bigr)f_j(X_i) + \int\Biggl( \sum_{j=1}^{M}
\bigl(\bar{\lambda}^{(2)}_j + \rho\eta_j\bigr)f_j(x) \Biggr)^2 \,dx \\
&&{} + 8rL\sum_{j=1}^{M}\bigl| \bar{\lambda}^{(2)}_j + \rho\eta_j\bigr|
=C.
\end{eqnarray*}
By taking the derivative with respect to $\rho$ of $F(\rho)$, we
obtain that, for all $ 0 <\rho< 1$,
\begin{eqnarray*}
F^{\prime}(\rho) &=&
-\frac{2}{n} \sum_{i=1}^{n}\sum_{j=1}^{M}\eta_jf_j(X_i) +
8rL\sum_{j=1}^{M}\eta_j \operatorname{sign}\bigl( \bar{\lambda}^{(2)}_j +
\rho
\eta_j\bigr)
\\
&&{} + 2{\int} \Biggl( \sum_{j=1}^{M} \bigl(\bar{\lambda}^{(2)}_j +
\rho\eta_j\bigr)f_j(x) \Biggr) \Biggl(\sum_{j=1}^{M}\eta_jf_j(x)
\Biggr)\,dx = 0 .
\end{eqnarray*}
By continuity of $\rho\mapsto\bar{\lambda}^{(2)}_j + \rho\eta_j$,
there exists an open interval in $(0,1)$ on which
$\rho\mapsto\operatorname{sign}( \bar{\lambda}^{(2)}_j + \rho\eta_j)$ is
constant for all $j$. Therefore, on that interval,
\[
F^{\prime}(\rho) = 2\rho{\int} \Biggl(\sum_{j=1}^{M}\eta_jf_j(x)
\Biggr)^2\,dx + C',
\]
where $C'$ does not depend on $\rho$. This is compatible with
$F^{\prime}(\rho)=0, \forall0 <\rho< 1$, 
only if
\[
\sum_{j=1}^{M}\eta_jf_j(x) = 0\qquad \mbox{for all } x
\]
and, therefore,
\[
\sum_{j=1}^{M}\eta_j \langle f_j, f_k \rangle= 0\qquad \mbox{for all } k
\in\{1, \ldots, M\},
\]
which is the desired result. This completes the proof of the lemma.
\end{pf}
\end{appendix}

\printaddresses

\end{document}